\documentclass[11pt]{article}
\usepackage{latexsym}
\usepackage[dvips]{epsfig}
\newcommand{\binomial}[2]{\left(\!\!\!\begin{array}{c}
        {#1} \\
        {#2} 
\end{array}\!\!\!\right)}

\newcommand{\C}{\bf C}

\newcommand{\A}{{\mathcal A}}

\newcommand{\HH}{{\mathcal H}} 
\newcommand{\K}{{\mathcal K}} 
\newcommand{\Kbul}{{\mathcal K}^{\bullet}}
\newcommand{\Kbulntilde}{\widetilde{\K}_{n}^{\bullet}}
\newcommand{\Kb}{{\overline{\mathcal K}}}
\newcommand{\Ke}{{\K^{\bullet _1}}}

\newcommand{\Kk}{{\K^{\bullet _2}}}
\newcommand{\Km}{{\K^{\bullet _m}}}
\newcommand{\Ki}{{\K^{\bullet _i}}}

\newcommand{\Kun}{{\K^{\bullet _1}}}
\newcommand{\Kdeux}{{\K^{\bullet _2}}}
\newcommand{\Ktrois}{{\K^{\bullet _3}}}

\newcommand{\Kiw}{{\K_w^{\bullet _i}}}
\newcommand{\Kp}{{\K^{\Diamond}}}
\newcommand{\Kpntilde}{{\widetilde{\K}_{n}^{\Diamond}}}
\newcommand{\Kptilde}{{\widetilde{\K}^{\Diamond}}}
\newcommand{\Kpw}{{\K_w^{\Diamond}}}
\newcommand{\Ah}{ {\mathcal A}_3}
\newcommand{\hai}{\widehat{ {\mathcal A}}_{i}}
\newcommand{\Ae}{ {\mathcal A}_1}
\newcommand{\Ai}{ {\mathcal A}_i}
\newcommand{\Ar}{ {\mathcal A}_{i,{\bf r}}}
\newcommand{\Ahar}{ \widehat{{\mathcal A}}_{i,{\bf r}}}
\newcommand{\Arh}{\widehat{{\mathcal A}}_{i,{\bf r}}^3}
\newcommand{\Ark}{ \widehat{{\mathcal A}}_{i,{\bf r}}^2}

\newcommand{\Aj}{ {\mathcal A}_j}
\newcommand{\Ak}{ {\mathcal A}_2}
\newcommand{\Am}{ {\mathcal A}_m}

\newcommand{\mca}{$m$-ary cactus }
\newcommand{\mci}{$m$-ary cacti }
\newcommand{\mcim}{$m$-ary cacti}

\newcommand{\qed}{\mbox{$\Box$}\vspace{\baselineskip}}

\newtheorem{theorem}{Theorem}

\newtheorem{lemma}[theorem]{Lemma}
\newtheorem{proposition}[theorem]{Proposition}

\newtheorem{corollary}[theorem]{Corollary}

\newtheorem{observat}[theorem]{Remark}
\newenvironment{remark}{\begin{observat}\rm}{\end{observat}}

\newenvironment{proof}{\noindent {\bf Proof:}}{\hfill{{\qed}}}

\newcommand{\vanish}[1]{}
\if@titlepage
  \newenvironment{resume}{%
      \titlepage
      \null\vfil
      \@beginparpenalty\@lowpenalty
      \begin{center}%
        \bfseries R\'esum\'e
        \@endparpenalty\@M
      \end{center}}%
     {\par\vfil\null\endtitlepage}
\else
  
\fi

\input{psfig}
 
\begin{document}
\title{Enumeration of $m$-ary cacti}

\author{Mikl\'os B\'ona\footnote{Present address: Institute for Advanced Study, Princeton.}\\
        Michel Bousquet\\
        Gilbert Labelle\\
        Pierre Leroux\\
                     \\
             LACIM\footnote{With the partial support of FCAR (Qu\'ebec) and CRSNG (Canada).}\\
        Universit\'e du Qu\'ebec \`a Montr\'eal \\
        Montr\'eal, Qu\'ebec, H3C 3P8\\
        Canada
        }
\date{January 18, 1999}
\maketitle

\begin{abstract}

The purpose of this paper is to enumerate various classes of cyclically 
colored $m$-gonal plane cacti, called $m$-ary cacti.
This 
combinatorial problem is motivated by the topological classification of 
complex polynomials having at most $m$ critical values, studied by Zvonkin and 
others. We obtain
explicit formulae for both labelled and unlabelled 
$m$-ary cacti, according 
to i) the number of polygons, ii) the vertex-color distribution, iii) the
vertex-degree distribution of each color. We also enumerate  
\mci according to the order of their automorphism group.
Using a generalization of Otter's formula, we express 
the species of $m$-ary cacti in terms of rooted and of pointed 
cacti. A variant of the $m$-dimensional Lagrange inversion is then used
to enumerate these structures. The method of Liskovets for the enumeration of 
unrooted planar maps can also be adapted to $m$-ary cacti.
\end{abstract}
\section{Introduction}
A {\em cactus} is a connected simple graph in which each edge lies in exactly 
one elementary cycle. It is equivalent to say that all blocks (2-connected 
components) of a cactus are edges or elementary cycles, i.e., polygons. An 
{\em $m$-gonal cactus} ($m$-cactus for short) is a cactus all of whose 
polygons are $m$-gons, for some 
fixed $m\geq 2$. By convention, a 2-cactus is simply a tree. These graphs were 
previously called ``Husimi trees'', and their definition was given by Harary and 
Uhlenbeck \cite{hu1} following a paper 
by Husimi \cite{h16} on 
the cluster integrals in the theory of condensation in statistical mechanics. 
See also Riddell \cite{r4} and Uhlenbeck and Ford \cite{ufi}. Their 
enumeration according to the number of polygons was carried out in 
\cite{hu1}. See also Harary and Palmer \cite{harary} and \cite{h16}.

\begin{figure}[h]
 \begin{center}
  \epsfig{file=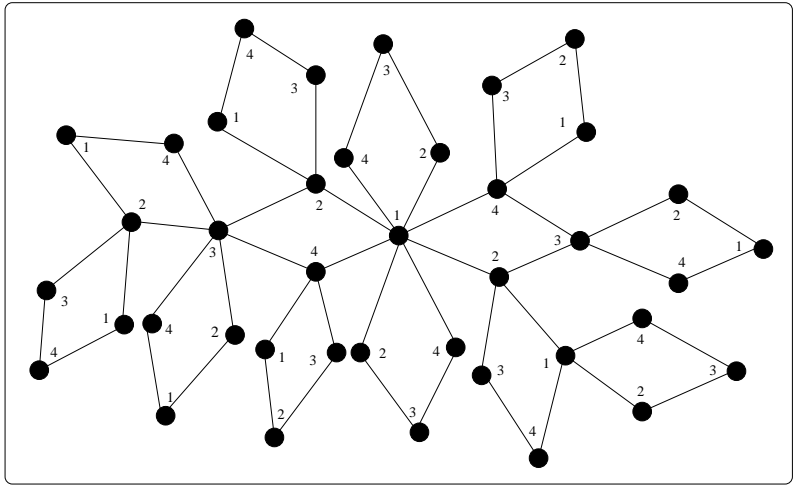, height=6.4cm}
  \caption{A quaternary cactus.}
  \label{4ary}
 \end{center}
\end{figure}

A {\em plane} $m$-cactus is an embedding of an $m$-cactus into the plane so 
that every edge is incident with the unbounded region. An $m$-{\em ary cactus}
is a plane $m$-cactus whose vertices are cyclically $m$-colored 
$1,2,\cdots, m$ counterclockwise within each $m$-gon. For technical reasons, 
we also consider a single vertex colored in any one of the $m$ colors to be 
an $m$-ary cactus. A quaternary ($m=4$) cactus is shown on Figure \ref{4ary}.

We define the {\em degree} of a vertex in a \mca to be the number of 
$m$-gons adjacent to that vertex. Note that it is twice the number of edges 
adjacent to the given vertex, for $m\geq 3$. Given an $m$-ary cactus $\kappa$, 
let $n_{ij}$ denote 
the number of vertices of color $i$ and degree $j$ of $\kappa$ and set 
$ {\bf n}_{i}=(n_{i0},n_{i1},n_{i2},\ldots )$. 
The {\em vertex-degree distribution} of $\kappa$ is given by the $m\times\infty$ 
matrix $N=(n_{ij})$, where $1\leq i\leq m$ and $j\geq 0$. Note that
 $n_i=\sum_jn_{ij}$ 
is the number of vertices of color $i$ and $n=\sum_in_i$ is the total number 
of vertices of $\kappa$. The \emph{vertex-color distribution} of $\kappa$ is defined 
to be the vector 
$\vec{\rm n}=(n_{1},n_{2},\ldots,n_{m})$. Also, let $p$ denote the number 
of polygons in $\kappa$.

For the quaternary cactus of Figure \ref{4ary}, the distributions are 
\[
{\rm {\bf n}}_1=(0,7,1,0,1,0,\cdots)=1^72^14^1,\;\; 
{\bf n}_2=(0,7,3,0,0,0,\cdots)=1^72^3,
\]
\[  
{\bf n}_3=(0,8,1,1,0,0,\cdots)=1^82^13^1,\;\;
{\bf n}_4=(0,9,2,0,0,0,\cdots)=1^92^2,
\]
\[
n_1=9,\; n_2=10,\; n_3=10,\; n_4=11,\; n=40,\; {\rm  and } \; p=13,\; 
\] 

Clearly, for any \mca with $n$ vertices and $p$ polygons we have
\begin{equation} \label{elso}  
        \sum_j jn_{ij} = p, \: \: \: \:
        \mbox{for all $i$,}
\end{equation}
since each polygon contains exactly one vertex of color $i$, and also 
\begin{equation} \label{masodik} 
        n=(m-1)p+1,     
\end{equation}
as one sees easily by induction on $p$.  

The goal of this paper is to enumerate various classes of \mci according to 
the number $n$ of vertices or $p$ of polygons, to the vertex color 
distribution $\vec{n}=(n_{1},n_{2},\ldots, n_{m})$, and to the vertex-degree 
distribution $N=(n_{ij})_{1\leq i\leq m,\ j\geq 0}$. 
The species we enumerate are the following:

\begin{enumerate}
\item $\K$,   the class of all \mcim, 
\item $\Ki$,  the class of \mci \emph{pointed} at vertex of color $i$ (see Figure \ref{coupure}), 
\item $\Kp$,  the class of \emph{rooted} (i.e., pointed at a polygon) cacti 
(see Figure \ref{polygon}),
\item $\A_{i}$, the class of \mcim, \emph{planted} at a vertex $v$ 
of color $i$, i.e., pointed at $v$ with a pair of half edges attached to $v$ 
        contributing to its degree (see Figure \ref{planted1}), 
\item $\Kb$,   the class of \emph{asymmetric} \mcim, 
\item $\K_{=s}$ and $\K_{\geq s}$,  the classes of \mci whose automorphism group is 
          of order $s$, and a multiple of $s$, respectively, where $s\geq 2$.
\end{enumerate}

\begin{figure}[h]
 \begin{center}
  \epsfig{file=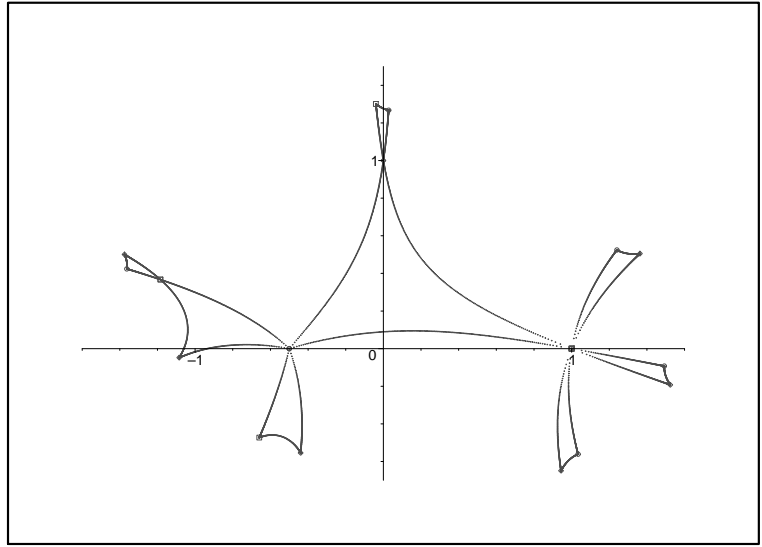, height=11cm, angle=270}
  \caption{Cactus associated to a polynomial 
           of degree 8, having three critical values.}
  \label{cactusmaple}
 \end{center}
\end{figure}

The motivation for the enumeration of \mci comes from the topological 
classification of polynomials having at most $m$ critical values. More 
precisely, two complex polynomials $p_{1}(z)$ and $p_{2}(z)$ are said to be 
\emph{topologically equivalent} if there exists two 
oriention-preserving homeomorphisms of the plane, $h_{1}$ and $h_{2}$, 
such that $h_{1}(p_{1}(z))=p_{2}(h_{2}(z))$. Also, a complex number $v$ is 
called a \emph{critical value} of the polynomial $p(z)$ if the equation 
$p(z)=v$ has at least one multiple root; all the roots of the equation are then called 
\emph{critical points}. Now if a polynomial $p(z)$ has $k$ critical values 
$\{v_{1}, v_{2}, \ldots, v_{k}\}$, we adjoin $m-k$ ``phoney'' critical 
values $\{v_{k+1}, \ldots, v_{m}\}$ and form a simple curvilinear $m$-gon 
joining the $m$ critical values $\{v_{1},\ldots, v_{m}\}$. Then the preimage 
under $p$ of this polygon yields an \mca whose vertex-degree distribution 
corresponds to the multiplicities of the critical points.
For example, Figure \ref{cactusmaple} shows the cactus corresponding to a degree $8$ 
polynomial $p(z)=c_{0}+c_{1}z+\ldots+c_{8}z^{8}$ having three critical 
values $v_{1}, v_{2}, v_{3}$, whose derivative is of the form 
$p'(z)=(z-b)(z-1)^{3}(z+\frac{1}{2})^{2}(z-i)$, where 
$b\in \C$ is chosen so that 
$p(b)=p(1)=v_{1}$, and $p(-\frac{1}{2})=v_{2}$, $p(i)=v_{3}$. 


This is a crucial 
step in the topological classification but the equivalence classes of 
polynomials are in fact represented
 by the orbits of \mci under the action of the braid group. See 
\cite{emhzz} and \cite{lesdeuxrusses} for more details. The enumeration of 
these orbits is an open problem. 

This work extends to general $m\geq 2$, previous results of Labelle and 
Leroux \cite{trees} on bicolored plane trees. It also extends results of 
Goulden and Jackson \cite{gj} on the enumeration of rooted \mcim. 
They show 
that rooted \mci with $p$ polygons, having vertex-degree distribution 
$N=(n_{ij})$ 
are in one-to-one correspondance with decompositions of the circular 
permutation $(1,2,\ldots,p)$ as the product $g_{1} g_{2} \cdots g_{m}$ 
of $m$ permutations, where $g_{i}$ has cyclic type 
$(1^{n_{i1}}2^{n_{i2}}\cdots)$.

\begin{figure}[h]
 \begin{center}
  \epsfig{file=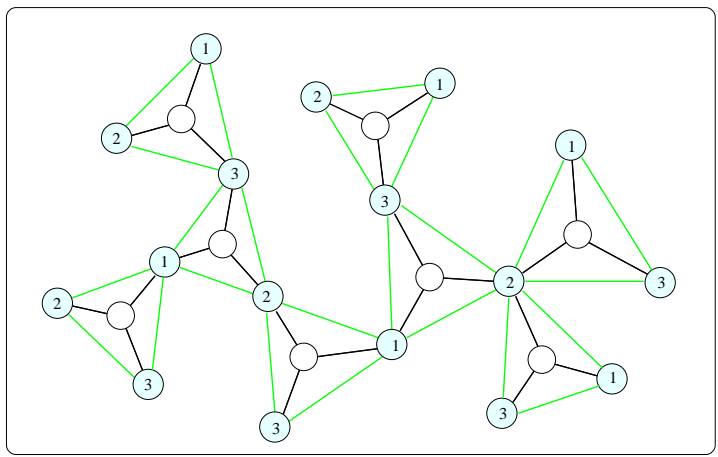, height=6cm}
  \caption{Tree-like structure of a ternary cactus}
  \label{cactus1bis}
 \end{center}
\end{figure}

In section 2, we state the main functional equations relating the 
various species of \mcim. We show that all these species can be expressed in 
terms of \emph{planted} \mci which, themselves, satisfy functional 
equations opening the way to Lagrange inversion. Of particular 
importance is a Dissymmetry Theorem which relates 
(ordinary) \mci to pointed and rooted \mcim. This 
theorem is closely related to the dissimilarity characteristic 
theorem for trees, due to Otter and extended to cacti 
by Harary and Norman \cite{hn1}. The tree-like structure of a cactus can be 
emphasized by using an equivalent representation, where a white 
($=$ color $0$) 
vertex is placed within each polygon, 
and joined to the vertices of the polygon, 
after which the edges of the polygons can be erased. This gives a 
bijection between \mci having $p$ polygons and 
$(1+m)$-colored trees having $p$ vertices of color $0$, all of degree $m$. 
The bijection is illustrated in Figure \ref{cactus1bis} 
for a ternary $(m=3)$ cactus. 

In section 3 we establish a particular form of multidimensional 
Lagrange inversion, which is well adapted to the present situation. It 
extends the previously known two-dimensional case, in the spirit of 
Chottin's formulae \cite{chottina} \cite{chottinb}, and use the crucial
observation due to Goulden and Jackson \cite{gj} that a certain Jacobian
matrix reduces to a rank-1 matrix. We then use these results in section 4 to 
enumerate both labelled and unlabelled \mcim, including the special 
classes of planted, pointed and rooted \mcim, according to the number of 
vertices
(or of polygons), to their vertex-color distribution and their 
vertex-degree distribution. We also enumerate \mci according to the order 
of their automorphism group, including the asymmetric ones. 

An alternate method can be used for the enumeration of unlabelled \mcim. It is 
based on a paper of Liskovets \cite{liskovets} on the enumeration 
of non rooted planar maps which uses the concept of quotient of a labelled planar 
map under an automorphism. See Bousquet \cite{bousq2}.

In the last section, we present some related enumerative results, concerning 
labelled {\em free} \mci and unlabelled plane 
$m$-gonal cacti having $p$ polygons. 
We also state a closely related result due to Bousquet-M\'elou 
and Schaeffer \cite{bms} on rooted $m$-ary  constellations, having $p$ polygons. 

Three tables are given at the end of the paper, containing numerical results 
which illustrate some of the formulas.

We have used the species formulation as a helpful unifying framework in this 
paper. A basic reference for the theory of species is the book 
\cite{species}. However, the paper remains accessible to anyone with a 
knowledge of P\'olya theory applied to graphical enumeration (see 
\cite{harary}). 

We would like to thank Sacha Zvonkin, for introducing us to the 
problem of cactus enumeration, and Robert Cori and Gilles Schaeffer for 
useful discussions.

\section{Functional equations for \mci }

\subsection{Vertex-color distribution}

We consider the class $\K$ of \mci as an $m$-sort species. This means that
an \mca is seen as a structure constructed on an $m$-tuple of sets 
$(U_1,U_2,\cdots ,U_m)$, the elements of $U_i$ being the (labels for) vertices
of color $i$. Moreover, the relabeling bijections, and in particular, the
automorphisms of \mci are required to preserve the sorts of elements, i.e.
the colors. Although we are interested in the enumeration of unlabelled 
 cacti, it is easier to establish the functional equations by giving bijections
between labelled structures. If we ensure that these bijections are natural,
that is, that they commute with any relabeling, thus defining isomorphisms of
species, then the consequences for both the labelled (exponential) generating
function
\begin{equation} \label{expo} 
        \K(x_1,x_2,\ldots ,x_m)=\sum_{n_1,n_2,\ldots,n_m}
        |\K[n_1,n_2,\ldots,n_m]|\frac{x_1^{n_1}\cdots x_m^{n_m}}
        {n_1!\cdots n_m!}
\end{equation}
and the unlabelled (ordinary) generating function
\begin{equation} \label{exp}
        \widetilde{\K}(x_1,x_2,\ldots,x_m)=\sum_{n_1,n_2,\ldots,n_m}
        \widetilde{\K}(n_1,n_2,\ldots,n_m) x_1^{n_1}\cdots x_m^{n_m}
\end{equation}
are automatic. Here $\K[n_1,n_2,\ldots,n_m]$ denotes the set of \mci over the 
multiset $[n_1]+[n_2]+\ldots +[n_m]$, with $[n]=\{1,2,\ldots,n\}$, and 
$\widetilde{\K}(n_1,n_2,\ldots,n_m)$ denotes the number of unlabelled \mci
having $n_i$ vertices of color $i$, for $i=1,\ldots, m$.

Note that the plane embedding of an \mca $\kappa$ is completely characterized 
by the specification, for each vertex $v$ of $\kappa$, 
of a circular permutation
on the polygons adjacent to $v$. We now present functional equations related 
to the $m$-sort species
 $\Ai$, of $m$-ary cacti, {\em planted} at a vertex of color $i$,
 $\Ki$, of $m$-ary cacti, {\em pointed} at a vertex of color $i$,  
 $\Kp$, of {\em rooted}  \mcim.    

The following notations are used:
$X_i$ denotes the species of singletons of sort (or color) $i$, 
 $\cal{C}$ denotes the species of (non-empty) circular
permutations, $L$ denotes the species of lists (linear orders) and 
$\widehat{\A}_i:=\prod_{j\neq i} \A_{j}$ denotes the product of all $\Aj$ 
except $\A_i$.
\begin{proposition} 
We have the following isomorphisms of species, for $i=1,\ldots, m$:
\begin{equation} \label{e}
        \Ai=X_iL(\widehat{\A}_i),
        \vspace{-0.7mm}
\end{equation}
\begin{equation} \label{f} 
        \Ki = X_i (1+{\cal{C}}({\widehat{\A}}_{i})), 
\end{equation}
\begin{equation} \label{g}  
        \Kp=\Ae \Ak \cdots \Am .
\end{equation}
\end{proposition}

\begin{figure}[h]
 \begin{center}
  \epsfig{file=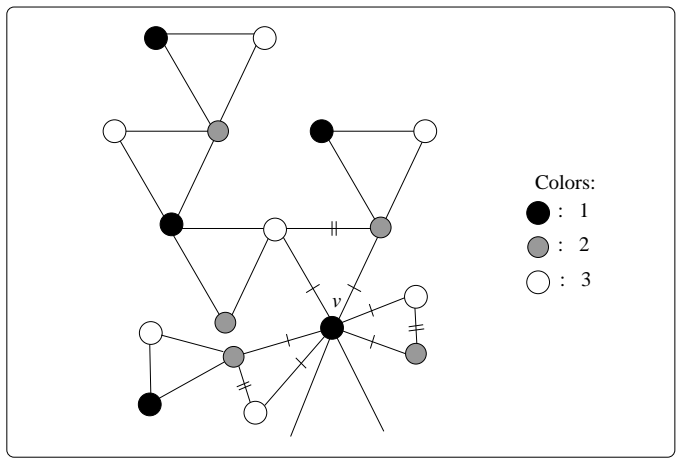,height=6cm}
  \caption{A planted ternary cactus.}
  \label{planted1}
 \end{center}
\end{figure}

\noindent{\bf Proof:}
The plane embedding of a planted \mca determines a
linear order on the neighboring polygons of the pointed vertex. If this 
vertex, say of color 1, is removed, each of these adjacent polygons can be 
simply decomposed into the product of  $m-1$ planted \mci with roots of color 
$2,3, \ldots, m$. Since this data completely specifies the planted cactus, we 
have equation (\ref{e}). See Figure \ref{planted1} for an illustration of 
the equation $\Ae=X_1L(\Ak \Ah )$ in the ternary case.

\begin{figure}[h]
 \begin{center}
  \epsfig{file=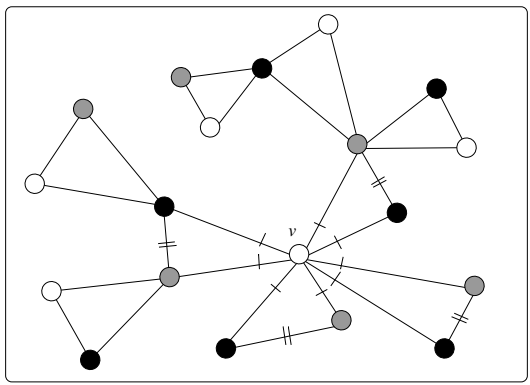, height=5cm}
  \caption{A ternary cactus pointed at vertex $v$.}
  \label{coupure}
 \end{center}
\end{figure}

\begin{figure}[h]
 \begin{center}
  \epsfig{file=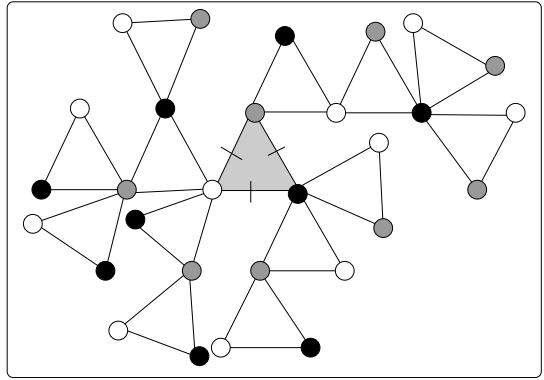, height=5cm}
  \caption{A rooted ternary cactus.}
  \label{polygon}
 \end{center}
\end{figure}

Equation (\ref{f}) is similar to (\ref{e}) except that for pointed cacti
the polygons adjacent to the pointed vertex can freely rotate around it. 
Figure \ref{coupure} illustrates the equation $\K^{\bullet _3}
=X_3(1+C(\Ae\Ak))$. Equation (\ref{g}) is immediate; see Figure \ref{polygon}.
\hfill $\Box$

Remark that equations (\ref{e}) and (\ref{g}) are essentially
due to Goulden and Jackson \cite{gj}. 

Recall that in a connected graph $g$, a vertex $x$ belongs to the 
{\em center} of $g$ if the maximal distance from $x$ to any other vertex is 
minimal. In particular, if $g$ is a cactus, then it is easy to see that the 
center of $g$ is either a single vertex or a polygon. Now let $\kappa$ be an 
$m$-ary cactus. In this case we define the {\em center} in a slightly 
different way: if the previous definition yields a vertex as the center of 
$\kappa$, 
then we 
leave this definition unchanged. 
If the previous definition yields a polygon $p$ 
as the center 
of $\kappa$, then we take the color-1 vertex of $p$ to be the center of 
$\kappa$. So now the center of an \mca  is always a vertex. 

\begin{theorem} \label{dissymth} {\sc Dissymmetry theorem for $m$-ary cacti}.
There is an isomorphism of species 
\begin{equation}
        \Ke + \Kk + \cdots + \Km = \K + (m-1) \Kp. \label{dissim}
\end{equation}
\end{theorem}
\begin{figure}[h]
 \begin{center}
  \epsfig{file=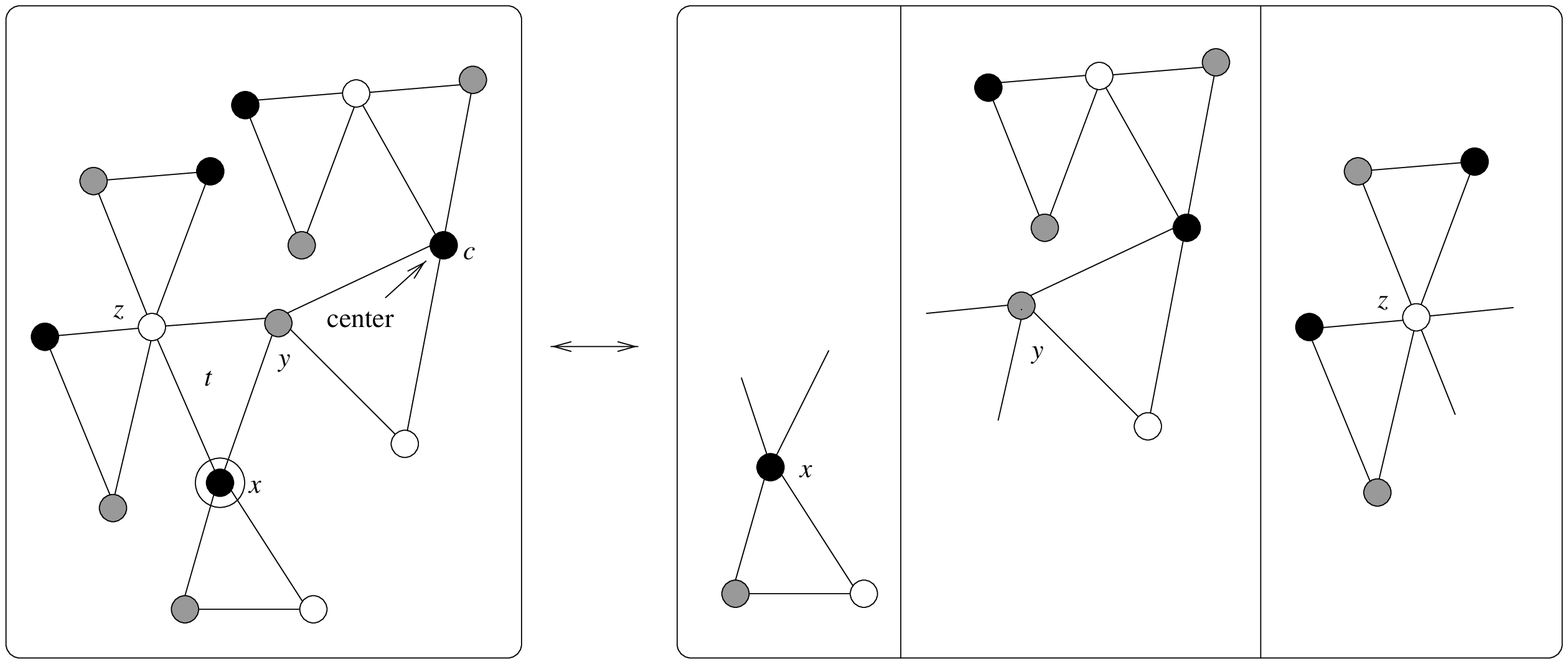, height=6cm}
  \caption{${\cal K}^{\bullet _1}+{\cal K}^{\bullet _2}
  +{\cal K}^{\bullet _3}={\cal K}+2{\cal K}^\Diamond$ $(m=3)$}
  \label{dissym2}
 \end{center}
\end{figure}
\begin{proof} For clarity, we prove the theorem for $m=3$, that is we 
establish an isomorphism 
$\Kun+\Kdeux+\Ktrois=\K+2\Kp$, the proof for
general $m$ being analogous.

The left hand side corresponds to  cacti which have been pointed at a
vertex, of color 1, 2, or 3. The first term of the right hand side 
corresponds to cacti which have been pointed in a canonical way, at their 
center. So what remains to construct is a natural bijection from triangular 
cacti pointed {\em not} in their center onto two cases of 
$\A_{1}\A_{2}\A_{3}$-structures.

Suppose that a ternary cactus $\kappa$ has been pointed at a vertex $x$ 
of color 1 which is 
different from the center $c$ of $\kappa$  (see Figure \ref{dissym2}). 
Let the shortest path from $x$ to 
$c$ start with the edge $e=\{x,y\}$, and let $t$ be the unique triangle
containing $e$. Then we cut the three edges of 
$t$ and thus separate the cactus into three smaller cacti  which are planted in a 
vertex of color 1, 2 and 3 respectively. We thus obtain an $\Ae \Ak \Ah $-structure. It is easy to see that we could have obtained this structure in another
way. Indeed, if the vertices of $t$ are $x,y$ and $z$, then pointing the 
cactus at $z$ would have given the same decomposition. So this operation 
does define a map into $2 \Ae \Ak \Ah  $.

To see that the algorithm is reversible, take any 3-tuple of cacti which are
planted in vertices $x$, $y$ and $z$ of color 1, 2, and 3 respectively. Join 
$x$, $y$ and $z$ by a triangle to get a cactus, and look for its center $c$. 
If $c$ comes from the component of $x$, then we can point either $y$ or 
$z$ in the cactus, if $c$ comes from the component of $y$, then we 
can point either 
$x$ or $z$ and finally, if $c$ comes from the component of $z$, then we 
can point either $x$ or $y$. It is then a simple matter to number each 
of these cases in order to make the correspondence bijective, completing the 
proof. \phantom{aaaaaaaa}
\end{proof}

\begin{corollary}
The species $\K$ of \mci can be written as 
\begin{eqnarray} \label{h} 
        \K&=&\sum_{i=1}^m\Ki-(m-1)\Kp \nonumber\\
          &=&\sum_{i=1}^m X_i(1+{\cal}{C}(\widehat{\A}_{i}))
               -(m-1)\prod_{i=1}^m\Ai.
\end{eqnarray}
\hfill $\Box$
\end{corollary}

The consequences for the labelled and unlabelled generating functions then
follow from general principles. For $i=1,\ldots, m$, we have, with
$\mathbf{x}=(x_{1},x_{2},\ldots,x_{m})$, 
\begin{equation} \label{lagen}
        \Ai({\bf x})=x_i\frac{1}{1-\hai({\bf x})},
\end{equation}

\begin{equation} \label{ungen}
        \widetilde{\Ai}({\bf x})=x_i\frac{1}{1-\widetilde{\hai}({\bf x})},
\end{equation}
from which it follows that $\widetilde{\Ai}({\bf x})=\Ai({\bf x})$ since 
$\widetilde{\widehat{{\cal A}}_i}=\widehat{\widetilde{{\cal A}}_i}$. This 
expresses the fact that planted cacti are asymmetric structures. Moreover,
\begin{equation}
        \Ki({\bf x})=x_i(1+\log \frac{1}{1-\hai({\bf x})}),
\end{equation}
\begin{equation} \label{29}
        \widetilde{\Ki}({\bf x})=x_i(1+\sum_{d\geq 1} 
        \frac{\phi(d)}{d}\log \frac{1}{1-\hai({\bf x}^d)}),
\end{equation}
where ${\bf x}^d:=(x_1^d,x_2^d,\ldots,x_m^d)$ and 
$\phi$ is the Euler function. 
We also have
\begin{equation} \label{hasz1}
        \Kp({\bf x})=\widetilde{\Kp}({\bf x})=\Ae({\bf x})\cdots \Am({\bf x})
\end{equation}
and finally, 
\begin{equation} \label{hasz2}
        \widetilde{\K}({\bf x})
        =\sum_{i=1}^m \widetilde{\Ki}({\bf x})-(m-1)\Kp ({\bf x}).
\end{equation}

\subsection{Vertex-degree distribution}

In order to enumerate \mci according to their degree distributions, we
 introduce weights in the form of monomials $w(\kappa)
=\prod_{i,j} r_{ij}^{n_{ij}}$
 with $i=1,\ldots, m$ and $j\geq 0$, for a cactus $\kappa$ having 
vertex-degree distribution
$N=(n_{ij})$. In other words, the variable 
$r_{ij}$ acts as a 
counter for (or marks) vertices of color $i$ and degree $j$. We also use the 
notation ${\bf r}_i$ to denote the sequence $(r_{i0},r_{i1},\ldots )$. We 
denote by $\K_w, \: \Kpw$, and $\Kiw$ the corresponding species of $m$-ary 
cacti, weighted in this manner. We denote by $\Ar$ the species of planted (at 
a vertex of color $i$) \mci similarly weighted by degree. The functional 
equations (\ref{e})--(\ref{h}) can then be extended as follows:
\begin{equation} \label{i}
        \Ar=X_i(r_{i,1}+r_{i,2}{\Ark}+r_{i,3}{\Arh}+\ldots )
\end{equation}
where $\widehat{{\mathcal A}}_{i,{\bf r}}=\prod_{j\neq i}\Ar,$
\begin{equation} \label{j}
        \Kiw=X_i(r_{i,0}+r_{i,1}{\cal{C}}_1(\Ahar)
+r_{i,2}{\cal{C}}_2(\Ahar)+\ldots ),
\end{equation}
where ${\cal{C}}_k$ denotes the species of circular permutations of length $k$, 
\begin{equation} \label{ka}
        \Kpw=\prod_{i=1}^m\Ar,
\end{equation}
and
\begin{equation} \label{k}
        \K_w=\sum_{i=1}^m \Kiw - (m-1) \Kpw.
\end{equation}

The important point here is that the weights behave multiplicatively, 
with respect to the operations of product and partitional composition.
The consequences for the labelled and unlabelled generating functions are
as follows:

\begin{equation} 
        \label{i1}
        \Ar({\bf x})=x_i(r_{i,1}+r_{i,2}\Ahar({\bf x})+
        r_{i,3}\Ahar^2({\bf x})+\ldots ), 
\end{equation}
\begin{equation}
        \label{j1}
        \widetilde{\Ar}({\bf x})=\Ar({\bf x}),
\end{equation}
\begin{equation}
        \label{k1}
        \Kiw({\bf x})
        =x_i(r_{i,0}+\sum_{h\geq 1}\frac{r_{i,h}}{h}\Ahar^h({\bf x})),
\end{equation}
\begin{equation}
        \label{k2}
        \widetilde{\Kiw}({\bf x})=x_i(r_{i,0}+\sum_{h\geq 1}\frac{r_{i,h}}{h}
        \sum_{d|h} \phi(d) \widehat{\A}_{i,{\bf r}^d}^{h/d} ({\bf x}^d) ),
\end{equation}
where ${\bf r}^d$ denotes the set of variables $\{r_{i,j}^d\}$, for 
$i=1,\ldots ,m$, $j\geq 0$. We also have 
\begin{equation} \label{k3}
        \Kpw({\bf x})
        =\widetilde{\Kpw}({\bf x})=\prod_{i=1}^m \Ar ({\bf x}), 
\end{equation}
and finally,
\begin{equation} \label{k4}
        \widetilde{\K_w}({\bf x})
        =\sum_{i=1}^m \widetilde{\Kiw}({\bf x})-(m-1)\Kpw({\bf x}).
\end{equation}

\subsection{One-sort \mci}
If neither the vertex-color nor the vertex-degree distribution are desired, but 
only the number of vertices or, equivalently, of polygons, then the enumeration 
is easier to carry out since one dimensional Lagrange inversion will suffice. 
Indeed, we can consider the various species of \mci introduced earlier as 
one-sort species, as Figure \ref{4ary} suggests. 
This means that the underlying 
set (of vertex labels) is independant of the colors and that the 
relabellings can be arbitrary, altough isomorphisms are still required 
to preserve 
colors. We use the same letters $\K,\; \Ki,\; \Kp,\; \Ai\;$ to denote 
these one-sort species. 
Equations (\ref{f})--(\ref{dissim}) are still valid in this context, 
with the following 
simplifications: first, all singleton species $X_{i}$ should be replaced 
by $X$; 
second, the addition of $1$ modulo $m$ to the colors induces isomorphisms 
of species 
$\A_{1}\cong \A_{2} \cong \ldots \cong \A_{m}$, and we write $\A$ for this 
common species, 
and also $\Kun \cong \Kdeux \cong \ldots \cong \Km$.

Equation (\ref{e}) then simplifies to 
\begin{equation}
\A=XL(\A^{m-1})=\frac{X}{1-\A^{m-1}}, \label{onzeune}
\end{equation}
which implies $\A = X+\A^{m}$. Moreover, equation (\ref{g}) reduces to
\begin{equation}
\Kp=\A^{m}=\A-X, \label{treizeune}
\end{equation}
while (\ref{f}) reduces to 
\begin{equation}
\Ki=X(1+C(\A^{m-1})). \label{douzeune}
\end{equation}

Finally, the dissymetry theorem for \mci takes the form

\begin{eqnarray}
\K & = & \K^{\bullet}-(m-1)\K^{\Diamond} \nonumber \\ 
   & = & m X (1+C(\A^{m-1}))-(m-1)(\A-X) , \label{quinzeune}
\end{eqnarray}  
where $\K^{\bullet}$ denotes the one-sort species of pointed at (any color) 
\mcim.

\section{Multidimensional Lagrange inversion techniques}

In this section we establish a special form of multidimensional Lagrange
inversion, which can be directly applied to $m$-ary cacti. First recall the
standard form, due to Good, (see Theorem 1.2.9, 1 of
\cite{gjlivre} or the equivalent formula (28b) of \cite{species}).

\begin{theorem} \label{lagrange} {\sc Good's Lagrange Inversion Formula.}
Let $A_1,A_2,\cdots, A_m$ be formal power series in the variables 
$x_1,x_2,\cdots,x_m$ such that the relations $A_i=x_iR_i(A_1,A_2,\cdots,A_m)$ 
are satisfied for all $i=1,\cdots,m$. Then for any formal power series 
$F(t_1,t_2,\cdots,t_m)$ we have:
\begin{equation} \label{elagrange}
        [x_1^{n_1}\cdots x_m^{n_m}]F( A_1(x),\cdots,A_m(x))
        =[t_1^{n_1}\cdots t_m^{n_m}]
        F({\bf t})|K({\bf t})|R_1^{n_1}({\bf t})\cdots
        R_m^{n_m}({\bf t}),
\end{equation}
where ${\bf t}=(t_1,t_2,\cdots,t_m)$ and $K({\bf t})$ is the $m\times m$ 
matrix whose $(i,j)$-th entry is
\begin{equation} \label{kdef} 
        K({\bf t})_{ij}
        =\delta_{ij}-\frac{t_j}{R_i({\bf t})}\cdot\frac{\partial 
        R_i({\bf t})}{\partial t_j}.  
\end{equation}
\hfill $\Box$
\end{theorem}

There is a particularly simple two-dimensional case of this formula,
the alternating case, which we call the {\em Chottin formula}. In the
papers \cite{chottina} \cite{chottinb}, Chottin worked extensively on the
two-dimensional Lagrange inversion and its combinatorial proof.

\begin{theorem} {\sc Chottin Formula}.
Let $A(x,y)$ and $B(x,y)$ be two formal power series satisfying the relations 
$A=x \Phi (B)$ and  $B=y \Psi (A)$, where $\Phi(t)$ and $\Psi(s)$ are
given formal power series. Then, for any non negative integers $\alpha$ and
$\beta$ we have:
\begin{equation} \label{chottin}
        [x^ny^m]A^{\alpha}B^{\beta}
        =(1-\frac{(n-\alpha)(m-\beta)}{nm})[s^{n-\alpha}t^{m-\beta}]
        \Phi^n(t)\Psi^m(s),\:\:\:\:\:  n\geq 1, \:\: m \geq 1. 
\end{equation}
\hfill $\Box$
\end{theorem} 

We extend this result into $m$ dimensions. 
\begin{theorem} {\sc Generalized Chottin formula}. \label{general}
Let $A_1, A_2,\ldots,A_m$ be formal power series in the variables 
$x_1,x_2,\ldots ,x_m$ such that for $i=1,\ldots,m$, the relations 
$A_i=x_i\Phi_i(\widehat{A}_i)$ are satisfied, where the $\Phi_i$ are given 
formal power series of one variable, and $\widehat{A}_i=\prod_{j\neq i}A_j$. 
Also let $n_1,\ldots ,n_m$  be integers $\geq 1$ and let $\alpha_1,\ldots,
\alpha_m$ be nonnegative integers. Set $n=\sum_{i=1}^m n_i$
and $\alpha=\sum_{i=1}^m \alpha_i$. Suppose that the following coherence
conditions are satisfied,
\[ 
        n_i\geq \alpha_i, \:\:\:\:\:\: \frac {n-\alpha}{m-1}=\beta  \:\:\:
        \mbox{is an integer},
\]
and set $\beta_i=\beta-n_i+\alpha_i$.
Then
\begin{equation} 
        [x_1^{n_1}\cdots x_m^{n_m}]A_1^{\alpha_1}\cdots A_m^{\alpha_m}
        =D\cdot[s_1^{\beta_1}\cdots s_m^{\beta_m}]\Phi_1^{n_1}(s_1)
        \cdots \Phi_m^{n_m}(s_m),
\end{equation}
where  
\begin{equation} \label{constante} 
        D=\prod_{i=1}^m (1+\frac{\beta_i}{n_i})-
        \sum_{j=1}^m \frac{\beta_j}{n_j} \prod_{i\neq j}
        (1+\frac{\beta_i}{n_i}).
\end{equation}

\end{theorem}

\begin{proof} We use Theorem \ref{lagrange} with $R_i(t_1,\ldots, t_m)=
\Phi_i(\widehat{t_i})$, where $\widehat{t}_i=\prod_{j\neq i} t_j$. We take 
advantage of some useful observations made by Goulden and Jackson in 
\cite{gj} to compute the determinant $|K({\bf t})|$. Indeed, for 
$i=1,2,\cdots , m$, we have
\[
        t_j\frac{\partial R_i}{\partial t_i}=0, 
\]  
as $R_i({\bf t})$ does not depend on $t_i$, and for $j\neq i$,  
\[
        t_j\frac{\partial R_i}{\partial t_j}
        =\widehat{t}_i\Phi'_i(\widehat{t}_i)
\] 
which is independent of $j$. We set $\Psi_i(\widehat{t}_i)
=\widehat{t}_i\Phi'_i(\widehat{t}_i)$ and write $\Psi_i=\Psi_i(\widehat{t}_i), 
\Phi_i=\Phi_i(\widehat{t}_i)$.

The definition of $K({\bf t})$ then yields, after routine transformations,
\[
        |K({\bf t})|
        =\frac{\prod_{i=1}^m (\Phi_i+\Psi_i)}{\prod_{i=1}^m\Phi_i} \cdot
        |\delta_{ij}-\frac{\Psi_i}{\Psi_i+\Phi_i}|.
\]
Let $M_{ij}=-\frac{\Psi_i}{\Psi_i+\Phi_i}$ and note that the rank of $M$ is
1 since all its columns are equal. So, by the Sherman-Morrison formula 
\cite{trace} we have $|I+M|=1+ \mbox{trace}(M)$. Therefore, the previous 
equation yields
\begin{equation} \label{matrix}
        |K({\bf t})|=\frac{\prod_{i=1}^m (\Phi_i+\Psi_i)}
        {\prod_{i=1}^m\Phi_i} \cdot
        (1-\sum_{i=1}^m \frac{\Psi_i}{\Psi_i+\Phi_i}).
\end{equation}
It follows from the Lagrange inversion formula (\ref{elagrange}) that 
\begin{eqnarray}
        [x_1^{n_1}\cdots x_m^{n_m}]A_1^{\alpha_1}
        \cdots A_m^{\alpha_m}\!\!\!
        &\!\!\!=\!\!\!&\!\![t_1^{n_1}\cdots t_m^{n_m}]t_1^{\alpha_1}
        \cdots t_m^{\alpha_m} \cdot |K({\bf t})|\prod_{i=1}^m
        \Phi_i^{n_i}\nonumber \\
        &=\!&\!\![t_1^{n_1-\alpha_1}
        \cdots t_m^{n_m-\alpha_m}] 
        (\prod_{i=1}^m \Phi_i^{n_i-1}(\Phi_i+\Psi_i))  
        (1-\sum_{i=1}^m \frac{\Psi_i}{\Psi_i+\Phi_i}).\label{longue}
\end{eqnarray}
Now let us define the coefficients $c_{i,\beta_i}$ by 
\begin{equation} \label{beta11} 
        \Phi_i^{n_i}(\widehat{t}_i)
        =\sum_{\beta_i\geq 0} c_{i,\beta_i}\widehat{t}_i^{\beta_i},
\end{equation}
which implies, by the definition of $\Psi_i$ that
\begin{equation} \label{beta12}
        \Phi_i^{n_i-1}(\widehat{t}_i)\Psi_i(\widehat{t}_i)
        =\sum_{\beta_i\geq 0}\frac{\beta_i}
        {n_i} c_{i,\beta_i} \widehat{t}_i^{\beta_i}.
\end{equation}
Recall that $n=\sum_{i=1}^m n_i$ and $\alpha=\sum_{i=1}^m\alpha_i$.  
Then $t_1^{n_1-\alpha_1}\cdots t_m^{n_m-\alpha_m}
=\widehat{t_{1}}^{\beta_1}\cdots\widehat{t_{m}}^{\beta_m}$ if and only if 
$\beta-\beta_i=n_i-\alpha_i$ for all $i$,where
$\beta=\sum_{i=1}^m\beta_i$. Summing these equations for 
$i=1,\ldots,m$ yields 
$(m-1)\beta=n-\alpha$ and $\beta_i=\beta-n_i+\alpha_i$.
We then conclude that (\ref{longue}) equals 
\[
        [\widehat{t}_1^{\beta_1}\cdots \widehat{t}_m^{\beta_m}]
        (\prod_{i=1}^m \Phi_i^{n_i-1}(\Phi_i+\Psi_i))  
        \cdot (1-\sum_{i=1}^m \frac{\Psi_i}{\Psi_i+\Phi_i})=
        D \cdot [s_1^{\beta_1}\cdots s_m^{\beta_m}]
        \Phi_1^{n_1}(s_1) \cdots \Phi_m^{n_m}(s_m),
\]
\noindent where $D$ is given by (\ref{constante}), completing the proof.
\end{proof}

The following special cases are particularly useful: 
\begin{enumerate}
\item 
        $\alpha_1=\alpha_2=\cdots =\alpha_m=1$, with the condition that 
        $(n-1)/(m-1)=p$ is a positive integer. Then we find that $\beta=p-1$, 
        $\beta_i=p-n_i$ and $D=p^{m-1}/\prod_{i=1}^m n_i$, and we have 
        \begin{equation} \label{chot1} 
                [x_1^{n_1}\cdots x_m^{n_m}]A_1\cdots A_m
                =\frac{p^{m-1}}{\prod_{i=1}^m n_i}\cdot[s_1^{p-n_1}\cdots 
                s_m^{p-n_m}]\Phi_1^{n_1}(s_1)\cdots 
                \Phi_m^{n_m}(s_m).
        \end{equation}
\item 
        $\alpha_1=0$, $\alpha_2=\cdots =\alpha_m=k\geq 1$, with the condition 
        that $(\sum_i a_i)/(m-1)=q$ is an integer. Then we find that 
        $\beta=q-k$, $\beta_1=q-a_1-k$, $\beta_i=q-a_i$, for $i=2,\cdots m$, 
        and that $D=q^{m-2}k/\prod_{i\neq 1} a_i$, and
        we have, writing $\widehat{A}_1=A_2A_3\cdots A_m$,
        \begin{equation} \label{kettes}
                [x_1^{a_1}\cdots x_m^{a_m}]\widehat{A}_1^k({\bf x}) 
                =\frac{q^{m-2}k}{\prod_{i\neq 1} a_i} \cdot 
                [s_1^{q-a_1-k}  s_2^{q-a_2} \cdots  s_m^{q-a_m}]
                \Phi_1^{a_1}(s_1) \cdots \Phi_m^{a_m}(s_m).
        \end{equation}
\item 
        Under the condition that $(\sum_i a_i)/(m-1)=q$ is an integer, it 
        follows from  (\ref{kettes}) that for any formal power series $F(s)$ 
        we have 
        \begin{equation} \label{c3}
                [x_1^{a_1}\cdots x_m^{a_m}]F(\widehat{A}_1)=\frac{q^{m-2}}
                {\prod_{i\neq 1} a_i}\cdot
                [s_1^{q-a_1-1}  s_2^{q-a_2} \cdots  s_m^{q-a_m}]F'(s_1)
                \Phi_1^{a_1}(s_1) \cdots \Phi_m^{a_m}(s_m).
        \end{equation}
\end{enumerate}
 
\section{Enumeration of \mci}

\subsection{Coherence conditions}
As observed in the introduction, there are some coherence conditions on the 
statistics of an $m$-ary cactus. 
We now state necessary and sufficient conditions for the 
existence of an $m$-ary cactus. 
The first one concerns the relationship between the number 
of vertices and the number of polygons. It is easily proved by induction on 
$p$. 

\begin{lemma} \label{lemmaun}
There exists an \mca having $n$ vertices and $p$ polygons if and only if 
\[
n=p(m-1)+1.
\]
\hfill{\qed}
\end{lemma}

\begin{lemma} \label{szin}
Let $\vec{n}=(n_1,n_2,\cdots ,n_m)$ be a vector of nonnegative integers and 
set $n=\sum_i n_i$.
There exists an \mca having $n$ vertices, $p$ polygons and vertex-color 
distribution $\vec{n}$ if and 
only if
\begin{enumerate}
        \item $p=(n-1)/(m-1)$ is an integer, 
        \item $p\geq 1$ $\Rightarrow$ $n_i \leq p$, for $i=1,\cdots,m$.
\end{enumerate}
\end{lemma}
\begin{proof}
The conditions are clearly necessary.
Sufficiency is proved by induction on $p$. If $p=0$, then $n=1$, and we
have a 1-vertex cactus. If $p\geq 1$, then all components of $\vec{n}$ are
strictly positive since otherwise, supposing for example that $n_1=0$, we find
\[
        n=\sum_{i=2}^m n_i \leq (m-1)p=n-1, 
\]
a contradiction. Hence we have $n_i\geq 1$, for all $i$. If $p=1$, then $n_i
=1$  for all $i$, and we have a cactus with a single polygon. If $p>1$, we
must have $n_i<p$ for some $i$, since otherwise $n=mp$ and $n=p(m-1)+1$ 
leads to a contradiction. Assume, say, $n_m<p$ and define a new vector
$\vec{n}'$ by $n_m'=n_m$ and $n_i'=n_i-1$ for $i=1,\cdots,m-1.$ This
vector $\vec{n}'$ satisfies the conditions 1 and 2 with $(n'-1)/(m-1)=
p-1$ and we can apply the induction hypothesis to construct a cactus with 
vertex distribution $\vec{n}'$. It suffices then to add a new polygon to this
cactus, attached to any existing vertex of color $m$ to obtain a cactus with 
vertex-color distribution equal to $\vec{n}$. 
\end{proof}

Observe that when conditions 1 and 2 are satisfied, 
$p\geq 1\Rightarrow n_i\geq 1$ for all $i$.

\begin{lemma}\label{lemmatrois}
Let $N=(n_{ij})_{1\leq i\leq m,j\geq 0}$ be an $m\times \infty$ matrix of 
non negative integers, and set $n=\sum_{ij}n_{ij}$. There exists an \mca having 
$n$ vertices and $p$ polygons and whose vertex-degree distribution is given 
by the matrix $N$ if and only if 
\begin{enumerate}
\item $p=(n-1)/(m-1)$ is an integer,
\item $\sum_{j} jn_{ij}=p$, for all $i$, 
\item $p\geq 1\Rightarrow n_{i0}=0$, for all $i$.
\end{enumerate}
\end{lemma}

\begin{proof}
These conditions are clearly necessary. Sufficiency is again proved by induction on 
$p$. If $p=0$, then $n=1$ and we have a one vertex cactus. If $p\geq 1$, then 
we can prove that for all $i$, except possibly one, $n_{i}\geq 1$. 
Indeed conditions $2$ and $3$ imply that $n_{i}=\sum_{j}n_{ij}\leq p$. 
Then, if $n_{i1}=0$ for some $i$, we have $n_{i}\leq p/2$. If this occurs for two or more 
values of $i$, then 
$n=\sum_{i}n_{i} \leq (m-1)p=n-1$, a contradiction. If $p=1$, then $n_{i1}=1$  for all 
$i$ and we have a one polygon cactus. If $p>1$, then either one $n_{i1}=0$, say 
$n_{m1}=0$, or all $n_{i1}$ are $\geq 1$. In the first case there must be some 
$j\geq 2$ with $n_{mj}\geq 1$; in the second case, there must exist some $i$, say $i=m$, 
and some $j\geq 2$, with $n_{mj}\geq 1$ since otherwise $n_{i}=n_{i1}=p$ for all 
$i$ and $n=\sum n_{i} =mp=(m-1)p+p$, a contradiction. In either case we set 
$n_{i1}'=n_{i1}-1$ for $i\not= m$, $n_{mj}'=n_{mj}-1$, $n_{m,j-1}'=n_{m,j-1}+1$ and 
$n_{ij}'=n_{ij}$ for other $i,j$. Then the matrix $N'=(n_{ij}')$ satisfies the 
conditions of the lemma with $p'=p-1$ and we can apply the induction hypothesis to 
construct a cactus with vertex-degree distribution $N'$. 
It remains then to add a new polygon to this cactus, attached to any existing 
vertex of color $m$ and degree 
$j-1$ to obtain a cactus with vertex-degree distribution $N$.
\end{proof}
\subsection{Rooted or labelled \mci}

As observed earlier, the species $\Kp$ of rooted \mci is asymmetric. It follows that 
labelled \mci and rooted \mci are closely related. For example, in the one-sort case, 
we have 
\begin{equation}
p\K_{n}=\K_{n}^{\Diamond}=n!\widetilde{\K}_{n}^{\Diamond} \label{pe}
\end{equation}
where $\K_{n}$ and $\K_{n}^{\Diamond}$ denote the number of \mci and rooted \mcim, 
respectively, having $n$ labelled vertices, and $\widetilde{\K}_{n}^{\Diamond}$ 
denotes the number of unlabelled \mci with $n$ vertices, and where $p$ is the number 
of polygons.

\begin{theorem}
Let $p$ be a positive integer and set $n=p(m-1)+1$.
Then the numbers $\widetilde{\K}_{n}^{\Diamond}$, of rooted (unlabelled) \mcim, and 
$\K_{n}$, of labelled \mcim, having $n$ vertices (and $p$ polygons), are given by
\begin{equation}
\Kpntilde=\frac{1}{n}{{mp}\choose{p}} \label{mr} 
\end{equation}
and
\begin{equation}
\K_{n}=\frac{(n-1)!}{p}{{mp}\choose{p}} .
\end{equation}
\end{theorem}

\begin{proof}
It follows from (\ref{onzeune}) and (\ref{treizeune}) that the one-sort species $\A$ and 
$\Kp$ of planted and rooted \mci respectively satisfy 
$\A(x)=x/(1-\A^{m-1}(x))$ and $\Kptilde(x)=\Kp(x)=\A(x)-x$. 
The result follows easily from Lagrange inversion since 
\begin{eqnarray*}
\Kpntilde &=& [x^{n}](\A(x)-x)\\
          &=& \frac{1}{n}[t^{n-1}](1-t^{m-1})^{-n}\\
          &=& \frac{1}{n}[t^{(m-1)p}](1-t^{m-1})^{-((m-1)p+1)} \\
          &=& \frac{1}{n}{{mp}\choose{p}}.
\end{eqnarray*}
The second result then follows from (\ref{pe}).
\end{proof}

\begin{remark}
Formula (\ref{mr}) also represents the number of (unlabelled) 
$m$-ary ordered rooted trees having $p$ internal vertices and $n$ leaves. A 
direct bijection can be given between rooted \mci and $m$-ary ordered rooted trees, 
which also explains the functional equation $\A=X+\A^{m}$. 
See \cite{bousq2} and \cite{bousqart}.
\end{remark}

Suppose now that a vector $\vec{n}=(n_1,n_2,\ldots ,n_m)$ satisfies the 
conditions of Lemma \ref{szin}. Let $\K_{\vec{n}}$ denote the number 
of \mci over the multiset of vertices $([n_1],[n_2],\ldots ,[n_m])$, that is, 
of labelled cacti with vertex-color distribution $\vec{n}$. Similarly let 
$\K_{\vec{n}}^{\Diamond}$ denote the number of labelled rooted \mci with vertex 
distribution $\vec{n}$. Then we have
\begin{equation} \label{deux} 
       p\cdot\K_{\vec{n}}= \K_{\vec{n}}^{\Diamond}=(\prod_{i=1}^m n_i!)
\widetilde{\K}_{\vec{n}}^{\Diamond},
\end{equation}
where $\widetilde{\K}_{\vec{n}}^{\Diamond}$ is the number of unlabelled rooted cacti with 
vertex-color distribution $\vec{n}$. 
\begin{theorem} \label{y}
Let $\vec{n}=(n_1,n_2,\ldots ,n_m)$ be a vector of nonnegative integers 
satisfying the coherence conditions of {\em Lemma \ref{szin}}, with 
$p\geq 1$. Then the number of unlabelled rooted \mci having vertex 
distribution $\vec{n}$ is given by 
\begin{equation} \label{trois}
        \widetilde{\K}_{\vec{n}}^{\Diamond}= \frac{1}{p}\prod_{i=1}^m {{p \choose n_i}}.
\end{equation}
\end{theorem}
\begin{proof}\
Recall that $\widetilde{\K}^{\Diamond}({\bf x})=\Kp({\bf x})=\Ae({\bf x})\Ak({\bf x})
\cdots\Am({\bf x})$ and that the $\Ai({\bf x})$ satisfy functional equation
(\ref{lagen}). Hence we can use the special case 1 of the Generalized Chottin 
formula, that is, formula (\ref{chot1}), with $\Phi_i(s)=L(s)=1/(1-s)$ for all $i$. Hence 
we find that
\begin{eqnarray*} 
        \widetilde{\K}_{\vec{n}}^{\Diamond}&=&[x_1^{n_1}\cdots x_m^{n_m}]
        \Ae({\bf x})\cdots\Am({\bf x})\\
        &=&\frac{p^{m-1}}{\prod_{i=1}^m n_i}
        \prod_{i=1}^m [s_i^{p-n_i}](1-s_i)^{-n_i}\\
        &=&\frac{p^{m-1}}{\prod_{i=1}^m n_i}\prod_{i=1}^m{{p-1 \choose n_i-1}},
\end{eqnarray*}
which implies (\ref{trois}).
\end{proof}

Putting together equations (\ref{deux})--(\ref{trois}) yields the following:

\begin{corollary}
If the conditions of {\rm Lemma \ref{szin}} are satisfied, the number of labelled
\mci with vertex-color distribution $\vec{n}=(n_1,n_2,\cdots ,n_m)$ is given by
\begin{equation}
        \K_{\vec{n}}=p^{m-2}\prod_{i=1}^m(p-n_i+1)^{<n_i-1>},
\end{equation}
where $x^{<k>}$ denotes the rising factorial $x(x+1)\cdots (x+k-1)$.
\end{corollary}

\begin{remark} This extends to general $m\geq 2$ 
the formula $n_1^{<n_2-1>}n_2^{<n_1-1>}$ 
for the number of labelled plane bicolored trees with vertex-color distribution
$(n_1,n_2)$ (see formula (2.7) of \cite{trees}).
\end{remark}  

To find the number $\K_{N}$ of labelled \mci having vertex-{\em degree} distribution 
$N=(n_{ij})$, where $i=1,\cdots,m$ and 
$j\geq 0$, a similar approach can be followed. As for the vertex-color distribution, 
we have
\begin{equation}
        p\cdot\K_{N}=n_1!\cdots n_m!
        \widetilde{\K}_{N}^{\Diamond},
\end{equation}
where $n_{i}=\sum_{j}n_{ij}$ and $\widetilde{\K}_{N}^{\Diamond}$ denotes the number of (unlabelled) 
rooted \mci having vertex-degree distribution $N$. Recall that 
${\bf n}_{i}=(n_{i0}, n_{i1}, n_{i2}, \ldots)$ is the degree distribution  for vertices 
of color $i$. The following result, due to Goulden and Jackson \cite{gj}, expresses 
the number $\widetilde{\K_{N}^{\Diamond}}$ in terms of the multinomial coefficients 
${{n_{i}}\choose{{\bf n}_{i}}}$. 

\begin{theorem}{\rm \cite{gj}}
Let $N=(n_{ij})_{1\leq i\leq m, j\geq 0}$ be an $m \times \infty$ matrix of non 
negative integers satisfying the coherence conditions of 
{\rm Lemma \ref{lemmatrois}}, 
with $n=\sum_{ij} n_{ij}$ and $p=(n-1)/(m-1)\geq 1$. Then the number of rooted 
\mci having $n_{ij}$ vertices of color $i$ and degree $j$, is given by
\begin{equation}
\widetilde{\K}_{N}^{\Diamond}=\frac{p^{m-1}}{\prod_{i=1}^{m} n_{i}} \prod_{i=1}^{m} 
{{n_{i}}\choose{{\bf n}_{i}}}. \label{troisr}
\end{equation}
\end{theorem}
\begin{proof}\  Recall that  $\Kpw({\bf x})=\widetilde{\K}_{w}^{\Diamond}({\bf x})=
\prod_{i=1}^m \Ar ({\bf x})$ and also recall equations (\ref{i1}). 
Again, we use the generalized Chottin formula (\ref{chot1}), with  
\begin{equation}\label{defphii}
        \Phi_i(s)=\Psi_{{\bf r}_i}(s):= r_{i1}+r_{i2}s+r_{i3}s^2+\cdots
\end{equation}  
Then we have 
\begin{equation} \label{preuveGJ}
        \widetilde{\K}_{N}^{\Diamond}
        =[\prod_{i,j}r_{ij}^{n_{ij}}][\prod_i x_i^{n_i}]\Kpw({\bf x})
        =\frac{p^{m-1}}{\prod_i n_i}[\prod_{i,j}r_{ij}^{n_{ij}}]  
        [\prod_i s_i^{p-n_i} ]\prod_{i=1}^m\Phi_i^{n_i}(s_i),
\end{equation}
which implies (\ref{troisr}).
\end{proof}
\begin{corollary} 
The number $\K_{N}$ of labelled \mci having vertex-degree 
distribution $N$, assuming that the conditions
of {\rm Lemma \ref{lemmatrois}} are satisfied, with $p\geq 1$, is 
given by 
\begin{equation} 
        \K_{N}
        =p^{m-2}\prod_{i=1}^m(n_i-1)!{{n_i\choose {\bf n}_{i} }}.
\end{equation}\hfill{{\qed}}
\end{corollary}

\begin{remark}
It is well-known that the number of ways to label an unlabelled structure 
$\kappa$ over an underlying multiset $[n_{1},n_{2},\ldots, n_{m}]$ is 
$n_{1}!n_{2}!\cdots n_{m}!/\!\!\mid\!\!\mathrm{Aut}(\kappa)\!\!\mid$, where 
$\mathrm{Aut}(\kappa)$
denotes the (color-preserving) automorphism group of $\kappa$. It follows that 
\begin{equation}
\K_{N} = \sum_{\kappa} \frac{n_{1}!n_{2}!\ldots n_{m}!}{\mid \!
\mathrm{Aut}(\kappa)\!\mid}
\end{equation}
where the sum is taken over all unlabelled \mci $\kappa$ with vertex-degree 
distribution $N$. It also follows that 
\begin{equation}
\sum_{\kappa} \frac{1}{\mid\!\mathrm{Aut}(\kappa)\!\mid} 
= \ \frac{1}{p} \widetilde{\K}_{N}^{\Diamond}.
\end{equation}
This formula can be used, as in \cite{emhzz}, 
to check that all unlabelled cacti with a given 
degree distribution have been found.
\end{remark}

\subsection{Pointed \mci (unlabelled)}
Recall that $\Kbul=\Kbul(X)$ denotes the one-sort species of \mci which are pointed at a 
vertex of any color. We have
\[
\widetilde{\K}^{\bullet}(x)=\widetilde{\K}^{\bullet_{1}}(x)
+\ldots+\widetilde{\K}^{\bullet_{m}}(x)=m\widetilde{\K}^{\bullet_{1}}(x).
\]
\begin{theorem}
Let $p$ be a positive integer and set $n=p(m-1)+1$. Then the number $\Kbulntilde$ of 
pointed \mci having $n$ vertices (and $p$ polygons), is given by
\begin{equation}
\Kbulntilde=\frac{1}{p} \sum_{d\mid p}\phi(d){{pm/d}\choose{p/d}}, \label{***}
\end{equation}
where $\phi$ is the Euler function.
\end{theorem}
\begin{proof}
We have $\Kbulntilde=m\widetilde{\K_{n}}^{\bullet _1}$ and 
$\Ke=X(1+C(\A^{m-1}))$. By Lagrange inversion, we find for $p\geq 1, n\geq m$, 
\begin{eqnarray*}
\widetilde{\K_{n}}^{\bullet _1}&=&[x^{n}](\widetilde{\K_{n}}^{\bullet _1}(x)-x)\\
   &=& [x^{n}]x\sum_{d\geq 1} \frac{\phi(d)}{d} \log{\frac{1}{1-\A^{m-1}(x^{d})}} \\
   &=& \sum_{d\mid n-1} \frac{\phi(d)}{d}[x^{\frac{n-1}{d}}]\log{\frac{1}{1-\A^{m-1}(x)}}\\
   &=& \sum_{d\mid n-1} \frac{\phi(d)}{n-1}(m-1) [t^{\frac{n-d-1}{d}}]t^{m-2}
       (1-t^{m-1})^{-\frac{n+d-1}{d}}\\
   &=&  \sum_{d\mid p} \frac{\phi(d)}{p} [t^{\frac{p-d}{d}}](1-t)^{-\frac{n+d-1}{d}}\\
   &=& \sum_{d\mid p} \frac{\phi(d)}{p}{{\frac{pm}{d}-1}\choose{\frac{p}{d}-1}}\\
   &=& \frac{1}{mp}\sum_{d\mid p} \phi(d){{pm/d}\choose{p/d}},
\end{eqnarray*}
which completes the proof.
\end{proof}

We now wish to compute the numbers $\widetilde{\K}_{\vec{n}}^{\bullet _i}$ and 
$\widetilde{\K}_{N}^{\bullet_{i}}$ of 
(unlabelled) \mci pointed at a vertex of color $i$, with vertex-color distribution 
$\vec{n}$ and vertex-degree distribution $N$, respectively. 
For symmetry reasons, it is sufficient 
to consider the case $i=1$ since we have
\begin{equation}
\widetilde{\K}_{\vec{n}}^{\bullet _i}=\widetilde{\K}_{\sigma^{i-1}\vec{n}}^{\bullet _1}
 \ \mbox{ and } \
\widetilde{\K}_{N}^{\bullet_{i}}=\widetilde{\K}_{\sigma^{i-1}N}^{\bullet_{1}} \label{ksigma}
\end{equation}
where $\sigma$ denotes a cyclic shift of the components of 
$\vec{n}$ or of the rows of $N$, i.e.
\begin{equation}
(\sigma\vec{n})_{i}=n_{i+1} \mbox{ and } (\sigma N)_{ij}=n_{i+1,j} ,
\end{equation}
the sum $i+1$ being taken modulo $m$. We introduce the following notations:
\begin{equation}
\vec{e}_{k}=(\delta_{ki}),\; i=1,\ldots,m,\;\;
\mathbf{e}_{h}=(\delta_{jh})_{j\geq 0},\;\;
E_{r,s}=(\delta_{ir}\cdot \delta_{js})_{1\leq i\leq m,\; j\geq 0}. 
\end{equation}

\begin{theorem}
Let $\vec{n}=(n_1, n_2, \ldots, n_m)$ be a vector of non negative integers 
satisfying the coherence conditions of {\rm Lemma \ref{szin}}, with 
$n=\sum_{i}n_{i}$ and $p=(n-1)/(m-1)\geq 1$. Then the number of \mci pointed 
at a vertex of color $1$ and having vertex-color distribution $\vec{n}$ is given 
by
\begin{equation}
\widetilde{\K}_{\vec{n}}^{\bullet _1}=
\frac{p-n_{1}+1}{p^{2}} \sum_{d} \phi(d) {{p/d}\choose{(n_{1}-1)/d}}
\prod_{i\not= 1} {{p/d}\choose{n_{i}/d}}, \label{deuxp}
\end{equation}
where the sum is taken over all $d$ such that $d$ divides $p$ and  all components of 
$\vec{n}-\vec{e_{1}}$.
\end{theorem}
\begin{proof}
Recall equation (\ref{29}), with $i=1$. In what follows we use the special
case 3 of the generalized Chottin formula, i.e. (\ref{c3}), with $F(s)
=\log\frac{1}{1-s}$ and $a_1=(n_1-1)/d$, $a_2=n_2/d$,$\ldots$, $a_m=n_m/d$, 
so that $q=p/d$. We find
\begin{eqnarray*} 
        \widetilde{\K}_{\vec{n}}^{\bullet _1}
        & = &[x_1^{n_1}\!\cdots x_m^{n_m}]
            (\widetilde{\K}^{\bullet_{i}}({\bf x})-x_1)\\
        & = &[\mathbf{x}^{\vec{n}}]x_1\sum_{d\geq 1}\frac{\phi(d)}{d}
        \log\frac{1}{1-\widehat{A}_1({\bf x}^d)}\\
        &= & [{\mathbf x}^{\vec{n}-\vec{e}_1}]
        \sum_{d\geq 1}\frac{\phi(d)}{d} 
        \log \frac{1}{1-\widehat{A}_1({\bf x}^d)}\\
        &= & \sum_{d|\vec{n}-\vec{e}_1}
        \frac{\phi(d)}{d}[{\bf x}^{(\vec{n}-\vec{e}_1)/d}]
        \log \frac{1}{1-\widehat{A}_1({\bf x})}\\
        &= & \sum_{d|(p,\vec{n}-\vec{e}_1)}\phi(d)\frac{p^{m-2}}
        {\prod_{i=2}^m n_i }\cdot [s_1^{\frac{p-n_1-d+1}{d}}]
        (\frac{1}{1-s_1})^{\frac{n_1+d-1}{d}}
        \prod_{i=2}^m[s_i^{\frac{q-n_i}{d}}](\frac{1}{1-s_i})^{n_i/d}\\
        &= & \sum_{d|(p,\vec{n}-\vec{e}_1)}\phi(d)\frac{p^{m-2}}
        {\prod_{i=2}^m n_i }{(p-d)/d \choose (n_1-1)/d} \prod_{i=2}^m
        {(p-d)/d \choose (n_i-d)/d}, 
\end{eqnarray*}
which is equivalent to (\ref{deuxp}). 
\end{proof}
\begin{theorem}
Let $N=(n_{ij})$ be an $m\times\infty$ matrix of non negative integers satisfying the 
coherence conditions of {\rm Lemma \ref{lemmatrois}}, with $n=\sum_{ij}n_{ij}$ and 
$p=(n-1)(m-1)\geq 1$. Then the number of \mci pointed at a vertex of color $1$ and 
having $n$ vertices of color $i$ and degree $j$, is given by 
\begin{equation}
\widetilde{\K}_{N}^{\bullet _1}=\frac{p^{m-2}}{\prod_{i\not= 1} n_{i}} 
\sum_{h,d} \phi(d)
{{(n_{1}-1)/d}\choose{(\mathbf{n}_{1}-\mathbf{e}_{h})/d}}
\prod_{i\not= 1} {{n_{i}/d}\choose{\mathbf{n}/d}}, \label{troisp}
\end{equation}
where the sum is taken over all ordered pairs $(h,d)$ such that $n_{1h}\not= 0$ and 
$d$ divides $h, p$ and all components of $\mathbf{n}_{1}-\mathbf{e}_{h}$ and of 
$\mathbf{n}_{i}$ with $i\geq 2$. 
\end{theorem}

\begin{proof}
 Recall that $n_{i}=\sum_{j}n_{ij}$ and $\mathbf{n}_{i}=(n_{ij})_{j\geq 0}$.
We 
will use the special case 2 of the generalized Chottin formula, i.e. 
(\ref{kettes}), with $\Phi_i(s)=\Psi_{{\mathbf r}_i^d}(s)$ (see
(\ref{defphii})), $k=h/d$, 
$a_1=(n_1-1)/d$, $a_2=n_2/d$, $\ldots$, $a_m=n_m/d$, for which $q=p/d$, 
and $A_i(\mathbf{x})={\mathcal A} _{i,{\bf r}^d}(\mathbf{x})$. 
Since $p\geq 1$, we have, 
using formula (\ref{k2}) 
with $i=1$, 
\begin{eqnarray}
        \widetilde{\K}_{N}^{\bullet _1}
        & = & [\prod_{ij} r_{ij}^{n_{ij}}][\prod_i x_i^{n_i}]
        (\widetilde{{\cal K}}_w^{\bullet_1}({\bf x})-r_{1,0}x_1)
        \nonumber\\
        & = & [{\bf r}^{N}][{\bf x^n}]\sum_{h} 
        \frac{x_1r_{1,h}}{h} \sum_{d|h}\phi(d)
        \widehat{{\mathcal A}} _{1,{\bf r}^d}^{h/d}({\bf x}^d)
        \nonumber\\
        & = & \sum_{h} \frac{1}{h} [{\bf r}^{N-E_{1,h}}]
        \sum_{d|(h,\vec{n}-\vec{e}_1)}\phi(d) [x_1^{(n_1-1)/d}x_2^{n_2/d}\! 
        \cdots x_m^{n_m/d}]
        \widehat{\mathcal A}_{1,{\bf r}^d}^{h/d}({\bf x})
        \nonumber\\
        & = & \frac{p^{m-2}}{\prod_{i\neq 1}n_i}
        \sum_{h,d}\phi(d)[{\bf r}^{N-E_{1,h}}]
        [s_1^{\frac{p-n_1-h+1}{d}}]\Psi_{{{\bf r}_1^d}}^{\frac{n_1-1}{d}}(s_1)
        \prod_{i\neq 1}[s_i^{\frac{p-n_i}{d}}]
        \Psi_{{\bf r}_i^d}^{\frac{n_i}{d}}(s_i)
        \nonumber\\
        & = & \frac{p^{m-2}}{\prod_{i\neq 1}n_i}\sum_{h,d}\phi(d)  
        [{\bf r}_1^{\frac{\mathbf{n}_{1}-\mathbf{e}_{h}}{d}}][s_1^{\frac{p-n_1-h+1}{d}}]
        \Psi_{{\bf r}_1^d}^{\frac{n_1-1}{d}}(s_1)
        \prod_{i\neq 1}[{\bf r}_i^{\frac{{\bf n}_i}{d}}][s_i^{\frac{p-n_i}{d}}]
        \Psi_{{\bf r}_i}^{\frac{n_i}{d}}(s_i)
        \nonumber\\
        & = & \frac{p^{m-2}}{\prod_{i\neq 1}n_i}
        \sum_{d,h}\phi(d)
        {{(n_1-1)/d}\choose ({\bf n}_1-\delta^h)/d}\prod_{i=2}^m
        {{n_i/d}\choose {\bf n}_i/d}, 
        \label{bousq2}
\end{eqnarray}
where the summation is taken over all ordered pairs $(h,d)$ such that
$n_{1,h} \not= 0$ and $d$ divides $h$, $p$, ${\mathbf n}_1
-\mathbf{e}_{n}$, and all ${\mathbf n}_i$ with $i\geq 2$.  
\end{proof}

\subsection{\mci (unlabelled)}
In order to enumerate unlabelled and unrooted \mcim, two methods can be used. The first one 
uses the dissymetry theorem for cacti (see Theorem \ref{dissymth}) 
which expresses the species of \mci in terms of pointed and of rooted cacti; see below. 
The second is Liskovets' method for the enumeration of unlabelled planar maps 
\cite{liskovets}. 
It uses the Cauchy-Frobenius theorem (alias Burnside's Lemma) and the 
concept of quotient of a planar map under an automorphism; see \cite{bousq2} and 
\cite{bousqart} for the application of Liskovet's method to the enumeration of \mcim.
\begin{theorem}
Let $p$ be a positive integer and set $n=p(m-1)+1$. Then the number 
$\widetilde{\K}_{n}$ of (unlabelled) \mci having $n$ vertices (and $p$ polygons), 
is given by
\begin{equation}
\widetilde{\K}_{n}=\frac{1}{p}\left( \frac{1}{n}{{mp}\choose{p}}+
\sum_{{d\mid p}\atop{d<p}} \phi(p/d){{md}\choose{d}}\right)\!, \label{unec}
\end{equation}
where $\phi$ is the Euler function.
\end{theorem}
\begin{proof}
Using the dissymmetry formula  (\ref{quinzeune}) for one-sort \mci, we find
\begin{equation}
\widetilde{\K_{n}}=\widetilde{\Kbul_{n}}-(m-1)\Kpntilde
\end{equation}
and the result follows from (\ref{mr}) and (\ref{***}).
\end{proof}

See Table $3$ for some numerical values of $\widetilde{\K}_{n}$.

\begin{theorem}
Let $\vec{n}=(n_{1}, n_{2}, \ldots, n_{m})$ be a vector of non negative integers 
satisfying the coherence conditions of {\rm Lemma \ref{szin}}, with 
$n=\sum_{i}n_{i}$ and $p=(n-1)/(m-1)\geq 1$. Then the number 
$\widetilde{\K}_{\vec{n}}$ of (unlabelled) \mci having vertex-color distribution $\vec{n}$ is 
given by 
\begin{equation}
\widetilde{\K}_{\vec{n}}=\frac{1}{p^{2}}\left(
\prod_{i=1}^{m} {{p}\choose{n_{i}}}+
\sum_{i,d}\phi(d)(p-n_{i}+1){{p/d}\choose{(n_{i}-1)/d}}
\prod_{j\not= i} {{p/d}\choose{n_{j}/d}}\right)\!, \label{deuxc}
\end{equation}
where the sum is taken over all pairs $(i,d)$ such that $1\leq i\leq m,\; d>1$, 
 $d$ divides $p$ and all components of $\vec{n}-\vec{e_{i}}$.
\end{theorem}
\begin{proof}
Using the dissymmetry formula (\ref{hasz2}), we have 
\begin{equation}
\widetilde{\K}_{\vec{n}}=\sum_{i=1}^{m} 
\widetilde{\K}_{\vec{n}}^{\bullet_{i}}-(m-1)\widetilde{\K}_{\vec{n}}^{\Diamond}.
\end{equation}
The result follows from (\ref{trois}), (\ref{ksigma}) and (\ref{deuxp}).
\end{proof}

See Table $2$ for some numerical values of $\widetilde{\K}_{\vec{n}}$.
\begin{theorem}
Let $N=(n_{ij})$ be an $m\times\infty$ matrix of non negative integers satisfying the 
coherence conditions of {\rm Lemma \ref{lemmatrois}}, with $n=\sum_{ij}n_{ij}$ and 
$p=(n-1)/(m-1)\geq 1$. Then the number $\widetilde{\K}_{N}$ of (unlabelled) 
\mci having $n_{ij}$ vertices of color $i$ and degree $j$, is given by
\begin{equation}
\widetilde{\K}_{N}=p^{m-2}\left(
\prod_{i=1}^{m} \frac{1}{n_{i}}{{n_{i}}\choose{\mathbf{n}_{i}}}+
\sum_{i,h,d} \frac{\phi(d)}{\prod_{\ell\not= i}n_{\ell}}
{{(n_{i}-1)/d}\choose{(\mathbf{n}_{i}-\mathbf{e}_{h})/d}}
\prod_{\ell\not= i} {{n_{\ell}/d}\choose{\mathbf{n}_{\ell}/d}}\right)\!, \label{troisc}
\end{equation}
where the sum is taken over all triplets $(i,h,d)$ such that 
$n_{ih}\not= 0,\; d>1,\;$ and $d$ divides $h, p$ and all entries of the matrix 
$N-E_{ih}$.
\end{theorem}
\begin{proof}
The dissymmetry formula (\ref{k}) gives 
\begin{equation}
\widetilde{\K}_{N}= \sum_{i=1}^{m} \widetilde{\K}_{N}^{\bullet_{i}}-(m-1)
\widetilde{\K}_{N}^{\Diamond}.
\end{equation}
The result follows from (\ref{troisr}), (\ref{ksigma}) and (\ref{troisp}).
\end{proof}

See Table $1$ for some numerical values of $\widetilde{\K}_{N}$.
\subsection{Unlabelled \mci according to their automorphisms}
We first consider asymmetric \mcim, that is, cacti whose automorphism 
group is reduced to the identity. 
Let $\overline{K}$ denote the species of asymmetric \mcim. 
We have already observed that the species 
$\Kp$ of rooted \mci is asymmetric i.e. that $\overline{\K}^{\Diamond}=\Kp$. 
The dissymetry 
formulas (\ref{quinzeune}) and (\ref{dissim}), yields, in the one-sort case,
\begin{equation}
\overline{\K}=\overline{\K}^{\bullet}-(m-1)\Kp,
\end{equation}
and in the $m$-sort case,
\begin{equation}
\overline{\K}=\sum_{i=1}^{m} \overline{\K}^{\bullet_{i}}-(m-1)\Kp.
\end{equation}
Since $\Ki=X(1+C(\widehat{\A_{i}}))$, the enumeration of (unlabelled) 
$\bar{\Ki}$-structures uses the \emph{asymmetry index series} $\Gamma_{C}$ of the 
species $C$ of circular permutations, instead of the 
cycle index series $Z_{C}$ for the enumeration of unlabelled cacti
(see \cite{gilbert}, \cite{species}),
where 
\begin{equation}
\Gamma_{C}(x_1,x_2,\ldots)=\sum_{d\geq 1} \frac{\mu(d)}{d}\log{\frac{1}{1-x_{d}}},
\end{equation}
compared to the cycle index series
\begin{equation}
Z_{C}(x_1,x_2,\ldots )= \sum_{d\geq 1} \frac{\phi(d)}{d} \log{\frac{1}{1-x_{d}}},\nonumber
\end{equation}
where $\mu$ is the M\"obius function.
It follows that the enumeration formulas for asymmetric \mci will be very similar 
to those of unlabelled cacti. In fact it suffices to replace $\phi$ by $\mu$ in the 
formulas of the previous section. Hence we have the following theorem.
\begin{theorem}
Assume that the coherence conditions of {\rm Lemmas \ref{lemmaun}}, 
{\rm \ref{szin}} and {\rm \ref{lemmatrois}}
are satisfied, with $p\geq 1$. Then the corresponding enumerative formulas for
(unlabelled) asymmetric \mci are as follows:
\begin{eqnarray}
\overline{\K}_{n} & = & \frac{1}{p}\left(\frac{1}{n}{{mp}\choose{p}}+
                   \sum_{{d\mid p}\atop{d<p}}\mu(p/d){{md}\choose{d}}\right)\!,\label{una}\\
\overline{\K}_{\vec{n}} & =& \frac{1}{p^{2}}\left( \prod_{i=1}^{m} {{p}\choose{n_{i}}}+ 
                        \sum_{i,d}\mu(d)(p-n_{i}+1)
                       {{p/d}\choose{(n_{i}-1)/d}}\prod_{j\not= i}
                      {{p/d}\choose{n_{j}/d}} \right)\!, \label{deuxa}\\
\overline{\K}_{N} &= & p^{m-2}\left(
\prod_{i=1}^{m} \frac{1}{n_{i}}{{n_{i}}\choose{\mathbf{n}_{i}}}+
\sum_{i,h,d} \frac{\mu(d)}{\prod_{\ell\not= i}n_{\ell}}
{{(n_{i}-1)/d}\choose{(\mathbf{n}_{i}-\mathbf{e}_{h})/d}}
\prod_{\ell\not= i} {{n_{\ell}/d}\choose{\mathbf{n}_{\ell}/d}}\right)\!, \label{troisa}
\end{eqnarray}
where the summation ranges of {\rm (\ref{deuxa})} and {\rm (\ref{troisa})} are the same as for 
{\rm (\ref{deuxc})} and {\rm (\ref{troisc})}.
\hfill\qed
\end{theorem}

\begin{figure}[h]
 \begin{center}
  \epsfig{file=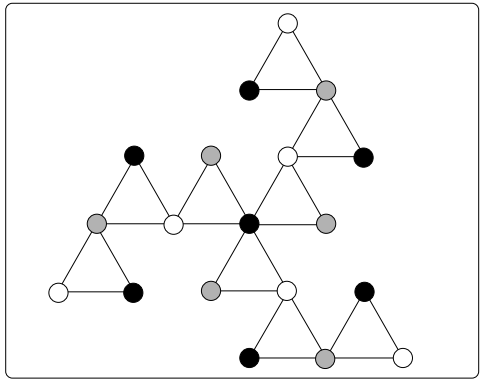, height=5cm}
  \caption{A ternary cactus with a symmetry of order 3.}
  \label{autom3}
 \end{center}
\end{figure}

We now consider \mci admitting at least one non trivial automorphism.
Since automorphisms are required to preserve colors, the only possibilities 
are rotations around a central vertex. See Figure 
\ref{autom3}. Observe that the order of such an automorphism must divide 
the number $p$ of polygons. Let $s\geq 2$ be an integer. Let $\K_{=s},$ and $\K_{\geq s},$ 
denote the 
species of \mci whose automorphism groups 
(necessarily cyclic) are of order $s$, and a 
multiple of $s$, respectively. Then, following the notations 
of \cite{trees}, section 3, we have
\begin{equation}
        \K_{=s} =\sum_{i=1}^m X_i C_{=s} (\widehat{\A}_{i}),
\end{equation}
\begin{equation}
        \K_{\geq s} =\sum_{i=1}^m X_i C_{\geq s} (\widehat{\A}_{i}).
\end{equation}

We can determine the unlabelled generating series $\widetilde{\K}_{\geq s} ({\bf x}) $ and
$\widetilde{\K}_{= s} ({\bf x}) $ by formulas (3.2) and (3.3) of \cite{trees}, essentially
due to Stockmeyer. See \cite{species}, Exercise 4.4.16, and \cite{stock}.
Extracting coefficients in these series is similar to the computations
of subsection 4.3. We find the following.
\begin{theorem}
Let $s\geq 2$ be an integer and assume that the coherence conditions of {\rm Lemmas
 \ref{lemmaun}}, {\rm \ref{szin}} and {\rm \ref{lemmatrois}} are satisfied, with $p$ a 
multiple of $s$. The corresponding enumerative formulas for (unlabelled) \mci 
whose automorphism groups are of order $s$, and a multiple of $s$, respectively, 
are as follows:

\begin{equation}
\widetilde{\K}_{=s,n}=\frac{s}{p}\sum_{d\mid \frac{p}{s}} \mu(d){{pm/sd}\choose{p/sd}}
\end{equation}
and
\begin{equation}
\widetilde{\K}_{\geq s,n}= \frac{s}{p} \sum_{d\mid \frac{p}{s}} 
\phi(d){{pm/sd}\choose{p/sd}};
\end{equation} 

\begin{equation}
\widetilde{\K}_{=s,\vec{n}}=\sum_{i=1}^{m} \frac{s(p-n_{i}+1)}{p^{2}}\sum_{d}\mu(d/s)
          {{p/d}\choose{(n_{i}-1)/d}}\prod_{j\not= i}{{p/d}\choose{n_{j}/d}},
\end{equation}
and
\begin{equation}
\widetilde{\K}_{\geq s,\vec{n}}= \sum_{i=1}^{m} \frac{s(p-n_{i}+1)}{p^{2}}\sum_{d}\phi(d/s)
          {{p/d}\choose{(n_{i}-1)/d}}\prod_{j\not= i}{{p/d}\choose{n_{j}/d}},
\end{equation}
the second summations being taken over all integers $d$ such that $s\!\!\mid\!\! d$ and 
$d \mbox{ divides } p$ and all components of  $\vec{n}-{\vec{e}}_{i}$;
\begin{equation}
\widetilde{\K}_{=s,N} = \sum_{i=1}^{m} \frac{p^{m-2}s}{\prod_{j\not= i} n_{j}}
\sum_{h,d}\mu(d/s) {{(n_{i}-1)/d}\choose{(\mathbf{n}_{i}-\mathbf{e}_{h})/d}}
\prod_{j\not= i} {{n_{j}/d}\choose{\mathbf{n}_{j}/d}},
\end{equation}
and 
\begin{equation}
\widetilde{\K}_{\geq s,N} = \sum_{i=1}^{m} \frac{p^{m-2}s}{\prod_{j\not= i} n_{j}}
\sum_{h,d}\phi(d/s) {{(n_{1}-1)/d}\choose{(\mathbf{n}_{i}-\mathbf{e}_{h})/d}}
\prod_{j\not= i} {{n_{j}/d}\choose{\mathbf{n}_{j}/d}},
\end{equation}
the second sommations  being taken over
all pairs of integers $h,d\geq 1$ such that $n_{ih}\not= 0, s\!\mid\! d$, and $d$ divides
$h$ and all entries in $N-E_{ih}$.
\hfill{{\qed}}
\end{theorem}

\section{Related enumerative results}
\subsection{Plane $m$-gonal cacti}
Let $\HH$ denote the one-sort species of plane $m$-gonal cacti (not $m$-colored).
The case of an isolated vertex is included. If $\HH^{\bullet}$ and $\A$ denote the species of 
pointed and of planted plane $m$-gonal cacti, respectively, then $\A$ coincides with the species introduced 
in section $2.3$, characterized by the functional equation $\A=XL(\A^{m-1})$, and 
$\HH^{\bullet}$ is isomorphic to the species $\K^{\bullet_{i}}$, for any $i$, that is, 
satisfies 
\begin{equation}
\HH^{\bullet}=X(1+\mathcal{C}(\A^{m-1})).
\end{equation}
See (\ref{onzeune}) et (\ref{douzeune}). However the species $\HH^{\Diamond}$ of rooted 
(at a polygon) plane $m$-gonal cacti is no longer asymmetric. In fact, we have 
\begin{equation}
\HH^{\Diamond} = \mathcal{C}_{m}(\A),
\end{equation}
where $\mathcal{C}_{m}$ denotes the species of circular permutations of length $m$. Another important 
difference resides in the form of the dissymetry theorem which is more closely 
related to that of free (non plane) $m$-gonal cacti. Indeed, we have (see \cite{hn1} and 
\cite{species}, 
(4.2.16) and Figure 4.2.5)
\begin{equation}
\HH^{\bullet}+\HH^{\Diamond}=\HH+\A\cdot\A^{m-1}
\end{equation}
from which we deduce, since $\A^{m}=\A-X$, that 
\begin{eqnarray}
\HH &=& \HH^{\bullet}+\HH^{\Diamond}-\A+X\nonumber\\
   &=& X(1+\mathcal{C}(\A^{m-1}))+\mathcal{C}_{m}(\A)-\A+X.
\end{eqnarray}
\begin{theorem}
Let $p$ be a positive integer and set $n=p(m-1)+1$. Then the numbers 
$\HH_{n}$ and $\widetilde{\HH}_{n}$ of labelled and unlabelled $m$-gonal cacti, 
repectively, having $n$ vertices (and $p$ polygons) are given by
\begin{equation}
\HH_{n}=\frac{(n-1)!}{mp}{{mp}\choose{p}},
\end{equation}
and
\begin{equation}
\widetilde{\HH}_{n}= \alpha_{n}+\beta_{n}-\gamma_{n},
\end{equation}
where
\begin{equation}
\alpha_{n}=\widetilde{\HH}_{n}^{\bullet} =
\frac{1}{mp} \sum_{d\mid p} \phi(\frac{p}{d}){{dm}\choose{d}},
\end{equation}
\begin{equation}
\beta_{n}=\widetilde{\HH}_{n}^{\Diamond} =
\frac{1}{mp} \sum_{d\mid (m,p-1)} \phi(d){{pm/d}\choose{(p-1)/d}},
\end{equation}
and
\begin{equation}
\gamma_{n}=\widetilde{\A}_{n}=\frac{1}{n} {{mp}\choose{p}}.
\end{equation}
\hfill\qed
\end{theorem}

In the case where $m=2$, we recover formulas of Walkup \cite{walkup} for the 
number of plane trees. See also Labelle and Leroux (\cite{trees}, (1.18)--(1.21)).
 It is also possible to derive similar formulas for the number of $m$-gonal 
plane cacti 
according to the vertex-degree distribution. See \cite{trees}, (1.23)--(1.26) where 
the computations have been carried out in the case $m=2$. Table $3$ contains numerical 
values of $\widetilde{\K}_{n}, \overline{\K}_{n}, \widetilde{\HH}_{n}$, for $n=(m-1)p+1$, 
and $m=2,\ldots, 7$. 

\subsection{Free (labelled) \mci}
A {\em free \mca} can be informally defined as an \mca without the plane
embedding. In other words, the $m$-gons attached to a vertex are free to take
any position with respect to each other. Denoting by ${\cal F}$ the species
of  free \mcim, we have the functional equations 
\begin{equation}
        {\cal F}^{\Diamond} = \Ae \Ak \cdots \Am,
\end{equation}
and, for $i=1,2,\cdots,m$,
\begin{equation} 
        \Ai=X_iE(\hai),
\end{equation}
where $E$ denotes the species of sets, for which 
\begin{equation}
        E(x)=e^x,\;
        \widetilde{E}(x)=(1-x)^{-1}\;\;{\rm and }\;\;\; Z_E(x_1,x_2,\cdots) 
        =\exp{(\sum_{i\geq 1}\frac{x_i}{i})}. \label{speciesE}
\end{equation}
The computations of subsection 4.1 for \emph{labelled} \mci according to
vertex-color distribution can be easily adapted to free \mcim. In particular, we
find the following result.
\begin{proposition} \label{ff}
Let $\vec{{\bf n}}=(n_1,\cdots,n_m)$ be a vector of positive integers 
satisfying the coherence conditions of {\rm Lemma \ref{szin}} with 
$p=(n-1)/(m-1)\geq 1$. Then the number 
${\cal F}(\vec{{\bf n}})$ of labelled free \mci having vertex-color 
distribution $\vec{{\bf n}}$ is given by 
\begin{equation} 
        {\cal F}(\vec{{\bf n}})
        =p^{m-2}\prod_{i=1}^m\frac{(n_i-1)!n_i^{p-n_i}}{(p-n_i)!}.
\end{equation}
\hfill\qed
\end{proposition}

This extends Scoins \cite{scoins} formula $n_1^{n_2-1}n_2^{n_1-1}$ for the number
of labelled bicolored free trees with vertex-color distribution $(n_1,n_2)$
to general $m\geq 2$. 

\subsection{Constellations}
%
%

\begin{figure}[h]
 \begin{center}
  \epsfig{file=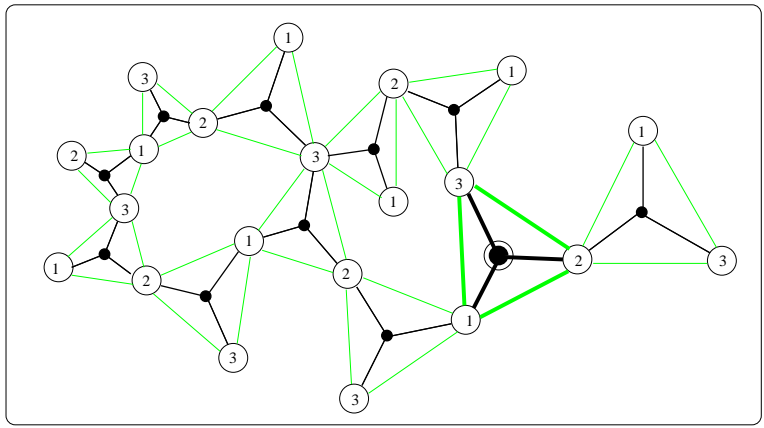}
  \caption{A rooted ternary constellation.}\newpage
  \label{constell}
 \end{center}
\end{figure}

Another combinatorial object closely related to $m$-ary cacti is an $m$-ary
constellation which is defined in a similar way as an $m$-ary cactus
 except that cycles of polygons are now allowed. Figure
\ref{constell} shows a typical ternary constellation which is rooted, 
that is, has a distinguished polygon. M. Bousquet-M\'elou and 
G. Schaeffer \cite{bms} have found that the number 
$\widetilde{C}^\Diamond(p)$ of unlabelled rooted $m$-ary constellations 
having $p$ polygons is given by
\begin{equation}
        \widetilde{C}^\Diamond(p)=
        \frac{(m+1)m^{p-1}}{((m-1)p+2)((m-1)p+1)}\binomial{mp}{p}.
\end{equation}

\section*{Tables}
\begin{center}
\begin{tabular}{|c|c|c|c|c|c|}\hline\hline
   $m$
   &\phantom{\LARGE{T}}$N$\phantom{\LARGE{T}} 
   & $\widetilde{\K}_{N}^{\bullet_{i}},\;i=1,\cdots,m$ 
   & $\widetilde{\K}_{N}^{\Diamond}$
   & $\widetilde{\K}_{N}$
   & $\overline{\K}_{N}$\\ \hline\hline
   2 &$(1^53^2,2^7)$               &(8, 7)           &14   &1    &1\\   \hline
   2 &$(1^22^24^1,1^22^4)$         &(76, 90)         &150  &16   &14\\  \hline
   3 &$(1^32^3,1^32^3,1^63^1)$     &(600, 600, 702)  &900  &102  &99\\  \hline
   3 &$(1^22^1,1^22^1,1^22^1)$     &(12, 12, 12)     &16   &4    &4\\   \hline
   3 &$(4, 1^{4},1^{4})$           &(1,1,1)          &1    &1    &0\\   \hline
   3 &$(2^{2},1^{2}2,1^{4})$        &(1,2,2)          &2    &1    &0\\   \hline
   3 &$(1^{1}3^{1},1^{2}2,1^{4})$     & (2,3,4)          &4    &1    &1 \\  \hline
   3 &$(1^22^2,1^22^2,1^42^1)$     &(54, 54, 69)     &81   &15   &12\\  \hline
   3 &$(1^32^14^1,1^32^3,1^72^1)$  &(600, 720, 960)  &1080 &120  &120\\ \hline
   3 &$(1^32^2,1^32^2,1^32^2)$     &(280, 280, 280)  &392  &56   &56\\  \hline
   3 &$(1^23^2,1^42^2,1^62^1)$     &(120, 180, 212)  &240  &32   &28\\  \hline
   3 &$(2^4,1^42^2,1^62^1)$        &(20, 30, 36)     &40   &6    &4\\   \hline
   3 &$(1^44^1,1^42^2,1^42^2)$     &(252, 300, 300)  &400  &52   &48\\  \hline
   3 &$(1^22^3,1^42^2,1^42^2)$     &(504, 600, 600)  &800  &104  &96\\  \hline
   4 &$(1^42^2,1^42^2,1^42^2,1^62^1)$ 
   &(6000, 6000, 6000, 7008) &8000 &1008 &992\\ \hline\hline
\end{tabular}
\end{center}
\begin{center}
        Table 1: The number of unlablelled \mci\\
        (rooted, plain, asymmetric) according to their vertex-degree distributions.
\end{center}


\begin{center}
\begin{tabular}{cc}
\begin{tabular}{|c|c|c|c|}\hline\hline
        \phantom{\LARGE{T}}$\vec{n}$\phantom{\LARGE{T}}    
        & $\widetilde{{\cal K}}^{\diamond}_{\vec{n}}$ 
        & $\widetilde{{\cal K}}_{\vec{n}}$
        & $\overline{{\cal K}}_{\vec{n}}$\\ \hline\hline
        $(7,7)$      & 226512       & 17424    & 17424\\      \hline
        $(5,6)$      & 5292         & 536      & 523\\        \hline
        $(6,6,7)$    & 28224        & 3138     & 3135\\       \hline 
        $(4,4,5)$    & 225          & 39       & 36\\         \hline
        $(5,6,8)$    & 10584        & 1176     & 1176\\       \hline
        $(5,5,5)$    & 1323         & 189      & 189\\        \hline
        $(4,6,7)$    & 1960         & 248      & 242\\        \hline
        $(5,6,6)$    & 5488         & 692      & 680\\        \hline
        $(3,4,4,5)$  & 50           & 10       & 10\\         \hline
        $(6,6,6,7)$  & 21952        & 2752     & 2736\\       \hline\hline
\end{tabular}&
\begin{tabular}{|c|c|c|c|}\hline\hline
        \phantom{\LARGE{T}}$\vec{n}$\phantom{\LARGE{T}}  
        & $\widetilde{{\cal K}}^{\diamond}_{\vec{n}}$ 
        & $\widetilde{{\cal K}}_{\vec{n}}$
        & $\overline{{\cal K}}_{\vec{n}}$\\ \hline\hline
        $(1,3,3)$    & 1            & 1       & 0\\    \hline
        $(2,2,3)$    & 3            & 1       & 1\\    \hline
        $(1,4,4)$    & 1            & 1       & 0\\    \hline
        $(2,3,4)$    & 6            & 2       & 1\\    \hline
        $(3,3,3)$    & 16           & 4       & 4\\    \hline
        $(3,3,5)$    & 20           & 4       & 4\\    \hline
        $(1,3,3,3)$  & 1            & 1       & 0\\    \hline
        $(2,2,3,3)$  & 3            & 1       & 1\\    \hline
        $(2,3,4,4)$  & 6            & 2       & 1\\    \hline
        $(4,4,4,4)$  & 125          & 25      & 25\\   \hline\hline
\end{tabular}
\end{tabular}
\end{center}
\begin{center}
        Table 2: The number of unlabelled \mci\\
        (rooted, plain, asymmetric) according  to their vertex-color distribution.
\end{center}

\begin{center}
\begin{tabular}{|c|c|c|c|c||c|c|c|c|}\hline
&\multicolumn{4}{|c||}{$m=2$}& \multicolumn{4}{|c|}{$m=3$}\\\hline\hline
$p$ & 
$\scriptstyle{n=p(m-1)+1}$&
$\widetilde{\K}_{n}$ & 
$\overline{\K}_{n}$ & 
$\widetilde{\HH}_{n}$& 
$\scriptstyle{n=p(m-1)+1}$&
$\widetilde{\K}_{n}$ & 
$\overline{\K}_{n}$ & 
\phantom{\LARGE{T}}$\widetilde{\HH}_{n}$\phantom{\LARGE{T}}\\ \hline
0 & 1 & 1   & 1  &  1  &1  &   1  &   1  &   1\\
1  &2  & 1  &  1  &  1 &3&  1  & 1 &  1\\
2  &3 &  2  &  0 &   1& 5 &  3  &  0   & 1\\
3  &4 &  3  &  1  &  2& 7& 6  &  3  &  2\\
4  &5 &  6   & 2  &  3& 9&19    &10  &  7\\
5   &6&  10  &  8  &  6&11 &57  &  54  &  19\\
6   &7 & 28  &  18   & 14&13 & 258  &  222  &  86\\
7  & 8&  63  &  61   & 34& 15&1110   & 1107   & 372\\
8   &9 & 190  &  170  &  95&17 & 5475  &  5346   & 1825\\
9   &10 & 546 &   538  &  280&19 &27429  &  27399  &  9143\\
10   &11 & 1708  &  1654 &   854&21 &143379  &  142770   & 47801\\
11   & 12& 5346  &  5344  &  2694&23 &764970   & 764967 &   254990\\
12   & 13& 17428  &  17252   & 8714&25 &4173906  &  4170672  &  1391302\\\hline\hline
\end{tabular}
\end{center}

\begin{center}
\begin{tabular}{|c|c|c|c|c||c|c|c|c|}\hline
&\multicolumn{4}{|c||}{$m=4$}& \multicolumn{4}{|c|}{$m=5$}\\\hline\hline
$p$ & 
$n$&
\phantom{\LARGE{T}}$\widetilde{\K}_{n}$\phantom{\LARGE{T}} & 
$\overline{\K}_{n}$ & 
$\widetilde{\HH}_{n}$& 
$n$&
$\widetilde{\K}_{n}$ & 
$\overline{\K}_{n}$ & 
$\widetilde{\HH}_{n}$\\ \hline
0 & 1&  1  &  1   & 1 &1 &  1  &  1  &  1\\
1  &4 &  1   & 1  &  1 &5 &  1  &  1  &  1\\
2  &7 &   4   & 0  &  1 &9 & 5  &  0  &  1\\
3   &10 & 10  &  6   & 3  &13&  15  &  10 &  3\\
4   & 13& 44    &28    &11  &17& 85   & 60   & 17\\
5   & 16& 197  &  193   & 52 &21 &  510   & 505  &  102\\
6   & 19& 1228  &  1140  &  307 &25 &  4051  &  3876   & 811\\
7   & 22& 7692  &  7688    &1936  &29& 33130  &  33125  &  6626\\
8   & 25& 52828   & 52364   & 13207  &33&  291925   & 290700  &  58385\\
9   & 28& 373636  &  373560   & 93496 &37 & 2661255   & 2661100  &  532251\\
10  & 31&  2735952  &  2732836  &  683988 &41 & 25059670  &  25049020 &   5011934\\
11  & 34&  20506258  &  20506254   & 5127163 &45 & 241724380  &  241724375 &   48344880\\
12   & 37& 156922676  &  156899748   & 39230669 &49 &  2379912355   & 2379812100  &  475982471\\
\hline\hline
\end{tabular}
\end{center}

\begin{center}
\begin{tabular}{|c|c|c|c|c||c|c|c|c|}\hline
&\multicolumn{4}{|c||}{$m=6$}& \multicolumn{4}{|c|}{$m=7$}\\\hline\hline
$p$ & 
$n$&
\phantom{\LARGE{T}}$\widetilde{\K}_{n}$\phantom{\LARGE{T}} & 
$\overline{\K}_{n}$ & 
$\widetilde{\HH}_{n}$& 
$n$&
$\widetilde{\K}_{n}$ & 
$\overline{\K}_{n}$ & 
$\widetilde{\HH}_{n}$\\ \hline
0 &1&  1 &  1 &  1&1& 1  &  1  &  1 \\
1 &6&  1 &  1 &  1&7& 1  &  1  &  1\\
2 &11&  6 &  0 &  1&13& 7   & 0  &  1\\
3 &16&  21 &  15 &  4&19& 28  &  21 &   4\\
4 &21&  146 &  110 &  25&25& 231  &  182  &  33\\
5 &26&  1101 &  1095  & 187&31& 2100  &  2093 &   300\\
6 &31&  10632 &  10326 &  1772&37& 23884  &  23394 &   3412\\
7 &36&  107062 &  107056  & 17880&43&  285390  &  285383  &  40770\\
8 &41&  1151802 &  1149126  & 191967&49&  3626295  &  3621150  &  518043\\
9  &46& 12845442 &  12845166  & 2141232&55& 47813815  &  47813367   & 6830545\\
10 &51&  147845706  & 147817170  & 24640989&61& 650367788   & 650302814  &  92909684\\\hline
\end{tabular}
\end{center}
\begin{center}
 Table 3: Number of $m$-ary and $m$-gonal cacti having $p$ polygons and $n$ vertices. 
\end{center}

\end{document}